\documentclass[12pt]{article}
\usepackage{verbatim,amsmath,amscd,amssymb,amsthm}
\begin{document}
\font\germ=eufm10
\def\ssl{\hbox{\germ sl}}
\def\slh{\widehat{\ssl_2}}
\title{
\bf
\Large Polyhedral Realizations of Crystal Bases for\\
Modified Quantum Algebras of Type $A$
}

\author{
H\textsc{oshino} Ayumu
\thanks
{
%supported in part by JSPS Grants in Aid for 
%Scientific Resarch
e-mail: a-hoshin@hoffman.cc.sophia.ac.jp
}
\\
Department of Mathematics, 
\\
Sophia University
\\
Tokyo 102-8554, Japan
\and
N\textsc{akashima} Toshiki
\thanks
{
%supported in part by JSPS Grants in Aid for 
%Scientific Resarch
e-mail: toshiki@mm.sophia.ac.jp
}
\\
Department of Mathematics, 
\\
Sophia University,
\\
Tokyo 102-8554, Japan
}
\date{}
\maketitle

\begin{abstract}
We describe the crystal bases of modified quantum algebras and
its connected component containing ``zero vector'' 
by the polyhedral realization method for the type $A_n$ and
$A^{(1)}_1$.
We also present the explicit form of the unique highest weight vector
in the connected component.
\end{abstract}

%\makeatletter
\def\aaa{@}
\vskip12pt
%\begin{tabular}{@{}c}

%\end{tabular}
\makeatother

\renewcommand{\labelenumi}{({\hbox{\roman{enumi}}})}
\font\germ=eufm10

\def\al{\alpha}
\def\beneme{\begin{enumerate}}
\def\beq{\begin{equation}}
\def\beqn{\begin{eqnarray}}
\def\beqnn{\begin{eqnarray*}}
\def\bigsl{{\hbox{\fontD \char'54}}}
\def\cd{\cdots}
\def\del{\delta}
\def\Del{\Delta}
\def\ei{e_i}
\def\eit{\tilde{e}_i}
\def\enepme{\end{enumerate}}
\def\ep{\epsilon}
\def\eeq{\end{equation}}
\def\eeqn{\end{eqnarray}}
\def\eeqnn{\end{eqnarray*}}
\def\fit{\tilde{f}_i}
\def\ft{\tilde{f}}
\def\ge{{\mathfrak g}}
\def\gl{\hbox{\germ gl}}
\def\hom{{\hbox{Hom}}}
\def\ify{\infty}
\def\io{\iota}
\def\kp{k^{(+)}}
\def\km{k^{(-)}}
\def\llra{\relbar\joinrel\relbar\joinrel\relbar\joinrel\rightarrow}
\def\lan{\langle}
\def\lar{\longrightarrow}
\def\lm{\lambda}
\def\Lm{\Lambda}
\def\mapright#1{\smash{\mathop{\longrightarrow}\limits^{#1}}}
\def\nd{\noindent}
\def\nn{\nonumber}
\def\ot{\otimes}
\def\op{\oplus}
\def\opi{\ovl\pi_{\lm}}
\def\ovl{\overline}
\def\plm{\Psi^{(\lm)}_{\io}}
\def\qq{\qquad}
\def\q{\quad}
\def\qed{\hfill\framebox[3mm]{}}
\def\QQ{\hbox{\bf Q}}
\def\qi{q_i}
\def\qii{q_i^{-1}}
\def\ran{\rangle}
\def\rlm{r_{\lm}}
\def\ssl{\hbox{\germ sl}}
\def\slh{\widehat{\ssl_2}}
\def\ti{t_i}
\def\tii{t_i^{-1}}
\def\til{\tilde}
\def\tt{{\hbox{\germ{t}}}}
\def\ttt{\hbox{\germ t}}
\def\uq{U_q(\ge)}
\def\uqm{U^-_q(\ge)}
\def\uqp{U^+_q(\ge)}
\def\uqmq{{U^-_q(\ge)}_{\bf Q}}
\def\uqpm{U^{\pm}_q(\ge)}
\def\uqmp{U^{\mp}_q(\ge)}
\def\uqq{U_{\bf Q}^-(\ge)}
\def\uqz{U^-_{\bf Z}(\ge)}
\def\util{\tilde{U}_q(\ge)}
\def\vep{\varepsilon}
\def\vp{\varphi}
\def\vpi{\varphi^{-1}}
\def\xii{\xi^{(i)}}
\def\Xiioi{\Xi_{\io}^{(i)}}
\def\wtil{\widetilde}
\def\what{\widehat}
\def\wpi{\widehat\pi_{\lm}}
\def\ZZ{\hbox{\bf Z}}
\newcommand{\fra}{\mathfrak a}
\newcommand{\frb}{\mathfrak b}
\newcommand{\frc}{\mathfrak c}
\newcommand{\frd}{\mathfrak d}
\newcommand{\fre}{\mathfrak e}
\newcommand{\frf}{\mathfrak f}
\newcommand{\frg}{\mathfrak g}
\newcommand{\frh}{\mathfrak h}
\newcommand{\fri}{\mathfrak i}
\newcommand{\frj}{\mathfrak j}
\newcommand{\frk}{\mathfrak k}
\newcommand{\frI}{\mathfrak I}
\newcommand{\fm}{\mathfrak m}
\newcommand{\frn}{\mathfrak n}
\newcommand{\frp}{\mathfrak p}
\newcommand{\fq}{\mathfrak q}
\newcommand{\frr}{\mathfrak r}
\newcommand{\frs}{\mathfrak s}
\newcommand{\frt}{\mathfrak t}
\newcommand{\fru}{\mathfrak u}
\newcommand{\frA}{\mathfrak A}
\newcommand{\frB}{\mathfrak B}
\newcommand{\frF}{\mathfrak F}
\newcommand{\frG}{\mathfrak G}
\newcommand{\frH}{\mathfrak H}
\newcommand{\frJ}{\mathfrak J}
\newcommand{\frN}{\mathfrak N}
\newcommand{\frP}{\mathfrak P}
\newcommand{\frT}{\mathfrak T}
\newcommand{\frU}{\mathfrak U}
\newcommand{\frV}{\mathfrak V}
\newcommand{\frX}{\mathfrak X}
\newcommand{\frY}{\mathfrak Y}
\newcommand{\frZ}{\mathfrak Z}
\newcommand{\rA}{\mathrm{A}}
\newcommand{\rC}{\mathrm{C}}
\newcommand{\rd}{\mathrm{d}}
\newcommand{\rB}{\mathrm{B}}
\newcommand{\rD}{\mathrm{D}}
\newcommand{\rE}{\mathrm{E}}
\newcommand{\rH}{\mathrm{H}}
\newcommand{\rK}{\mathrm{K}}
\newcommand{\rL}{\mathrm{L}}
\newcommand{\rM}{\mathrm{M}}
\newcommand{\rN}{\mathrm{N}}
\newcommand{\rR}{\mathrm{R}}
\newcommand{\rT}{\mathrm{T}}
\newcommand{\rZ}{\mathrm{Z}}
\newcommand{\bbA}{\mathbb A}
\newcommand{\bbC}{\mathbb C}
\newcommand{\bbG}{\mathbb G}
\newcommand{\bbF}{\mathbb F}
\newcommand{\bbH}{\mathbb H}
\newcommand{\bbP}{\mathbb P}
\newcommand{\bbN}{\mathbb N}
\newcommand{\bbQ}{\mathbb Q}
\newcommand{\bbR}{\mathbb R}
\newcommand{\bbV}{\mathbb V}
\newcommand{\bbZ}{\mathbb Z}

%%%%%%%%%%%%%%% Section 1 %%%%%%%%%%%%%%%%
\renewcommand{\thesection}{\arabic{section}}
\section{Introduction}
\setcounter{equation}{0}
\renewcommand{\theequation}{\thesection.\arabic{equation}}

Let 
$\uq:=<e_i,f_i,q^h>_{i\in I}$ ($I={1,2,\cd,n}$) 
be the quantum algebra associated with the 
symmetrizable Kac-Moody Lie algebra $\ge$ and 
$\uqm:= <f_i>_{i\in I}$ (resp. $\uqp:=<e_i>_{i\in I}$)
be the subalgebra of $\uq$. 
Kashiwara showed that 
the subalgebras $\uqpm$ has a unique 
crystal base $(L(\mp\ify),B(\mp\ify))$ and 
arbitrary integrable highest weight $\uq$-module 
$V(\lm)$ has a unique crystal baase 
$(L(\lm),B(\lm))$ (\cite{K1}). 

He also proved that a tensor product of crystal bases for
modules is again a crystal base for
tensor product of corresponding modules, which is
one of the most beautiful and usefull properties in theory of 
 crystal bases.

The term ``crystal '' implies a combinatorial notion 
abstracting the properties of
crystal base without assuming existence of 
the corresponding modules.
We shall see its examples in Sect.2.1.
We can define tensor product structure on crystals 
in a similar manner to crystal bases.
Indeed, some crystals (and their tensor products) are used 
to realize crystal base $B(\ify)$(\cite{NZ}) and
$B(\lm)$(\cite{N1}) by the polyhedral realization method.
Polyhedral realization of crystal bases is the method embedding
crystal bases in some infinite-dimensional vector space 
 and describing its image as a set of lattice points in certain convex 
polyhedron(\cite{N1},\cite{NZ}).

The modified quantum algebra 
$\util:=\bigoplus_{\lm \in P}\uq a_{\lm}$ 
(resp. $\uq a_{\lm}:=\uq/\sum\uq(q^h-q^{\lan h_i,\lm\ran}))$ 
is given by modifiying the Cartan part of $\uq.$
Lusztig showed that it has a 
crystal base $(L(\util),B(\util))
=(\oplus_{\lm}L(\uq a_{\lm}),\oplus_{\lm}B(\uq a_{\lm}))$
 (\cite{L}) and Kashiwara describe its important properties
(\cite{K4}). One of them is that 
the existence of 
the following isomorphism of crystals:
$$
B(\uq a_{\lm})\cong B(\ify)\ot T_{\lm}\ot B(-{\ify})
$$
where $T_{\lm}$ is the crystal given in Sect.2.2.
We have already the polyhedral realization of $B(\pm)$ and then 
we can get the polyhedral realization of $B(\uq a_{\lm})$ 
(Sect.4), which  is, in general, not connected, on the contrary, 
including infinitely many components. 
So,
next we try to describe some specific
connected component $B_0(\lm)$
containing $u_{\ify}\ot t_{\lm}\ot u_{-\ify}$
and the explicit
form of the unique highest weight vector in $B_0(\lm)$ by 
the  polyhedral realization method under certain assumptioin on the weight $\lm$
in the cases $\ge=A_n$ and $A^{(1)}_1$.

The organization of this paper is as follows:
in Sect.2., we review the theory of crystal bases and crytals.
We also prepare several ingredients to use in the subsequent 
sections, {\it e.g.,} crystals $B_i$, $T_\lm$ and explain 
polyhedral realization of $B(\pm\ify)$.
In Sect.3., we define the modified quantum algebra $\til\uq$ and 
see its crystal base and in Sect.4., its polyhedral realization 
is given. In Sect.5.,we consider the modified quantum 
algebra of type $A_n$.
We describe the highest weight vector
in the connected component $B_0(\lm)$ including 
$u_\ify\ot t_\lm\ot u_{-\ify}$ and the polyhedal realization
of $B_0(\lm)$ under certain condition on $\lm$.
In the last section, we treat the affine $A^{(1)}_1$-case.
In this case, we succeeded in presenting the explicit form 
of the connected component of $B_0(\lm)$ 
and the unique highest weight vector in $B_0(\lm)$ for a 
poisitive level weight $\lm$.

%%%%%%%%%%%%%%% Section 2 %%%%%%%%%%%%%%%%
\section{Crystal Base and Crystals}
\setcounter{equation}{0}
\renewcommand{\thesection}{\arabic{section}}
\renewcommand{\theequation}{\thesection.\arabic{equation}}

\theoremstyle{definition}
\newtheorem{def1}{Definition}[section]
\theoremstyle{plain}
\newtheorem{pro1}[def1]{Proposition}
\newtheorem{lem1}[def1]{Lemma}
\newtheorem{thm1}[def1]{Theorem}
\newtheorem{ex1}[def1]{Example}

%%%%%%%%%%%%%% 2.1. %%%%%%%%%%%%%
\subsection{Definition of Crystal Base and Crystal}

In this section, we shall review crystal bases and crystals.
We follow the notations and terminologies to \cite{N1}\cite{NZ}.

We fix a finite index set $I$. 
Let $A=(a_{ij})_{i,j\in I}$ be 
a generalized symmetrizable Cartan matrix,
$(\ttt,\{\al_i\}_{i\in I},\{h_i\}_{i\in I})$ be the associated 
Cartan data and $\ge$ be the associated Kac-Moody Lie algebra
where $\al_i$(resp. $h_i$) is called an simple root
(resp. simple coroot).
Let  $P$ be a weight lattice, $P^*$ be a dual lattice 
including $\{h_i\}_{\in I}$ and 
$Q:=\bigoplus_{i\in I}\QQ(q)\al_i$ be a root lattice.
Let $\uq:=\lan q^h,e_i,f_i\ran_{i\in I,h\in P^*}$ 
be the quantum algebra defined by the usual relations, which is an 
associative and Hopf algebra over the field $\QQ(q)$.
(We do not write down the Hopf algebra structure here.)

Now we give the definition of crystal base. Let ${\cal O}_{\rm int}$
be the category whose objects are $\uq$-module
that it has a weight space decomposition and 
for any element $u$, there exists positive integer 
$l$ such that $e_{i_1}\cd e_{i_k}u=0$ for any $k>l$.
It is well-known that the category ${\cal O}_{\rm int}$ is 
a semisimple category and all simple ojbects are parametrized 
by the set of dominant integral weights $P_+$. 
Let $M$ be a $\uq$-module  in $\cal O_{\rm int}$. 
For any $u\in M_\lm$ $(\lm\in P)$,
we have the unique expression:
\[
 u=\sum_{n\geq 0}f_i^{(n)}u_n,
\]
where $u_n\in {\rm Ker}\,e_i\cap M_{\lm+n\al_i}$. By using this,
we define the Kashiwara operators 
$\eit,\fit\in {\rm End}(M)$ $(i\in I)$,
\[
 \eit u:=\sum_{n\geq 1}f_i^{(n-1)}u_n,\,\,\qq
 \fit u:=\sum_{n\geq 0}f_i^{(n+1)}u_n.
\]
here note that we can define the Kashiwara operators 
$\eit,\fit\in {\rm End}(\uqpm)$ by the similar manner(\cite{K1}).
Let $A\subset \QQ(q)$ be the subring of rational functions which
are regular
 at $q=0$. Let $M$ be a $\uq$-module in $\cal O_{\rm int}$
(resp. $U^{\pm}_q(\ge)$).

\begin{def1}[\cite{K1}]
A pair $(L,B)$ is a crystal base of $M$ (resp. $\uqpm$), 
if it satisfies:
\begin{enumerate}
\item
$L$ is free $A$-submodule of $M$ (resp. $\uqpm$) and 
$M\cong \QQ(q)\ot_A L$ (resp. $\uqpm\cong \QQ(q)\ot_A L$).
\item
$B$ is a basis of the $\QQ$-vector space $L/qL$.
\item
$L=\oplus_{\lm \in P}L_{\lm}$, $B=\sqcup_{\lm \in P}B_{\lm}$ where
$L_{\lm} := L \cap M_{\lm}$, $B_{\lm} := B \cap L_{\lm}/qL$.
\item
$\eit L\subset L$ and
$\fit L\subset L$.
\item
$\eit B\subset B\sqcup\{0\}$ and
$\fit B\subset B\sqcup\{0\}$
(resp. $\fit B\subset B$) ( $\eit$ and $\fit$ acts on
$L/qL$ by (iv)).
\item
For $u,v\in B$,
$\fit u=v$ $\Longleftrightarrow$ $\eit v=u$.
\end{enumerate}
\end{def1}

The unit of subalgebra $\uqmp$ is denoted by
 $u_{\pm\ify}$. We set
\begin{eqnarray*}
&&L(\ify) := \sum_{i_j\in I,l\geq 0}
A\til f_{i_l}\cd \til f_{i_1}u_{\ify},\\
&&L(-\ify)  := \sum_{i_j\in I,l\geq 0}
A\til e_{i_l}\cd \til e_{i_1}u_{-\ify},
\\
&&B(\ify) :=
\{\til f_{i_l}\cd \til f_{i_1}u_{\ify}\,\,
{\rm mod}\,\,qL(\ify)\,|\,i_j\in I,l\geq 0\},\\
&&B(-\ify)  := 
\{\til e_{i_l}\cd \til e_{i_1}u_{-\ify}\,\,
{\rm mod}\,\,qL(-\ify)\,|\,i_j\in I,l\geq 0\}.
\end{eqnarray*}
\begin{thm1}[\cite{K1}]
A pair $(L(\pm\ify),B(\pm\ify))$ 
is a crystal base of $\uqmp$.
\end{thm1}

Now we introduce the notion {\it crystal}, which is obtained by
abstracting the combinatorial properties of crystal bases.
\begin{def1}
A {\it crystal} $B$ is a set endowed with the following maps:
\begin{eqnarray*}
&& wt:B\lar P,\\
&&\vep_i:B\lar\ZZ\sqcup\{-\infty\},\q
  \vp_i:B\lar\ZZ\sqcup\{-\infty\} \q{\hbox{for}}\q i\in I,\\
&&\eit:B\sqcup\{0\}\lar B\sqcup\{0\},
\q\fit:B\sqcup\{0\}\lar B\sqcup\{0\}\q{\hbox{for}}\q i\in I,\\
&&\eit(0)=\fit(0)=0.
\end{eqnarray*}
those maps satisfy the following axioms: for
 all $b,b_1,b_2 \in B$, we have
\begin{eqnarray*}
&&\vp_i(b)=\vep_i(b)+\lan h_i,wt(b)\ran,\\
&&wt(\eit b)=wt(b)+\al_i{\hbox{ if  }}\eit b\in B,\\
&&wt(\fit b)=wt(b)-\al_i{\hbox{ if  }}\fit b\in B,\\
&&\eit b_2=b_1 \Longleftrightarrow \fit b_1=b_2\,\,(\,b_1,b_2 \in B),\\
&&\vep_i(b)=-\ify
   \Longrightarrow \eit b=\fit b=0.
\end{eqnarray*}
\end{def1}
Indeed, if $(L,B)$ is a crystal base, then $B$ is a crystal.

\begin{def1}
\begin{enumerate}
\item
Let $B_1$ and $B_2$ be crystals. A {\it strict morphism} of crystals
$\psi:B_1\lar B_2$ is a map
$\psi:B_1\sqcup\{0\} \lar B_2\sqcup\{0\}$
satisfying the following: 
(1) $\psi(0)=0$. (2)If $b\in B_1$ and $\psi (b)\in B_2$, then
\[
\hspace{-30pt}wt(\psi(b))=wt(b),\q \vep_i(\psi(b))=\vep_i(b),\q
  \vp_i(\psi(b))=\vp_i(b).
\]
and the map $\psi$ commutes with all $\eit$ and $\fit$.
\item
An injective strict morphism is called an 
{\it embedding of crystals}. We
call $B_1$ is a subcrystal of $B_2$, if $B_1$ is a subset of $B_2$ and
becomes a crystal itself by restricting the data on it from $B_2$.
\end{enumerate}
\end{def1}

The following examples will play an important role
in the subsequent sections.
\begin{ex1}
\label{ex-tlm}
Let  $T_{\lm}:=\{t_{\lm}\}$ $(\lm\in P)$
be the crystal consisting of one element $t_\lm$ defined by
$wt(t_{\lm})=\lm,$
$\vep_i(t_{\lm})=\vp_i(t_{\lm})= -{\ify}$
, $\eit (t_{\lm})=\fit(t_{\lm})=0$.
\end{ex1}

\begin{ex1}
\label{ex-bi}
For $i\in I$, the crystal $B_i:=\{(x)_i\,: \, x \in \ZZ\}$ 
is defined by:
\beqnn
&& wt((x)_i)=x \al_i,\qq \vep_i((x)_i)=-x, \qq \vp_i((x)_i)=x,\\
&& \vep_j((x)_i)=-\ify,\qq \vp_j((x)_i)=-\ify \q {\rm for }\q j\ne i,\\
&& \til e_j (x)_i=\del_{i,j}(x+1)_i,\qq
\til f_j(x)_i=\del_{i,j}(x-1)_i.
\eeqnn
\end{ex1}
Note that as a set $B_i$ is identidied with the set of integers
$\bbZ$.

%%%%%%%%%%%%%%%%% 2.2 %%%%%%%%%%%%%
\subsection{Polyhedral Realization of $B(\pm\ify)$}
We review the polyhedral realization of the crystal $B(\pm\ify)$
following to \cite{NZ}.

We consider the following additive groups:
\[
\ZZ^{+\ify}
:=\{(\cd,x_k,\cd,x_2,x_1): x_k\in\ZZ
\,\,{\rm and}\,\,x_k=0\,\,{\rm for}\,\,k\gg 0\},
\]
\[
\ZZ^{-{\ify}}
:=\{(x_{-1},x_{-2}, \cd , x_{-k}, \cd): x_{-k} \in \ZZ 
\,\,{\rm and}\,\, x_{-k}=0 \,\,for \,\,k\gg 0 \}.
\]
We will denote by
 $\ZZ^{+{\ify}}_{\geq 0} \subset \ZZ^{+{\ify}}$
( resp. $\ZZ^{-{\ify}}_{\leq 0} \subset \ZZ^{-{\ify}}$ ) 
the semigroup of non-negative (resp. non-positive) integer
sequences. 
Take an infinite sequence of indices
$\io^+=(\cd,i_k,\cd,i_2,i_1)$ 
(resp. $\io^-=(i_{-1},i_{-2},\cd , i_{-k},\cd)$) from $I$ such that
\begin{equation}
{\hbox{
$i_k\ne i_{k+1}$ for any $k$, and $\sharp\{k>0 \,\,
{\rm (resp.\,\, k<0)}\,\,: i_k=i\}=\ify$ for any $i\in I$.}}
\label{seq-con}
\end{equation}
We can associate to 
 $\io^+$ ( resp. $\io^-$) a crystal structure on $\ZZ^{+{\ify}}$ 
(resp. $\ZZ^{-{\ify}}$) (see \cite{NZ}) and denote it by
 $\ZZ^{+{\ify}}_{\io^+}$ 
(resp. $\ZZ^{-{\ify}}_{\io^-}$). 
Let $B_i$ be the crystal given in Example\ref{ex-bi}. We 
obtain the following embeddings(\cite{K2}):
\begin{eqnarray*}
\Psi_{i}^+ : B(\infty)&\hookrightarrow &B(\infty)\ot B_i
\q
(u_{\ify}\mapsto u_{\ify}\ot (0)_i),\\
\Psi_{i}^- : B(-{\infty})&\hookrightarrow &B_i \ot B(-{\infty})
\q
(u_{-{\ify}}\mapsto (0)_i \ot u_{-{\ify}}).
\end{eqnarray*}
Iterating $\Psi_{i}^+$ (resp. $\Psi_{i}^-$)
according to $\io^+$ (resp. $\io^-$), we get 
the {\it Kashiwara embedding}(\cite{K2});
\begin{eqnarray}
&&\Psi_{\io^+}:B(\ify)\hookrightarrow \ZZ^{+{\ify}}_{\geq 0}
\subset \ZZ^{+{\ify}}_{\io^+}
\q
(u_{{\ify}}\mapsto (\cd,0,\cd,0,0,0)),
\label{kas+}\\
&&\Psi_{\io^-}:B(-{\ify})\hookrightarrow \ZZ^{-{\ify}}_{\leq 0}
\subset \ZZ^{-{\ify}}_{\io^-}
\q
(u_{-{\ify}}\mapsto (0,0,0,\cd,0,\cd)).
\label{kas-}
\end{eqnarray}

We consider the following infinite dimensional vector spaces
and their dual spaces:
\begin{eqnarray*}
\QQ^{+{\ify}}&:=&\QQ\ot_{\ZZ} \ZZ^{+\ify}=\{\vec{x}=
(\cd,x_k,\cd,x_2,x_1): x_k \in \QQ\,\,{\rm and }\,\,
x_k = 0\,\,{\rm for}\,\, k \gg 0\},\\
\QQ^{-{\ify}}&:=&\QQ\ot_{\ZZ} \ZZ^{-\ify}=\{\vec{x}=
(x_{-1},x_{-2},\cd,x_{-k},\cd): x_{-k} \in \QQ\,\,{\rm and }\,\,
x_{-k} = 0\,\,{\rm for}\,\, k \gg 0\},\\
 (\QQ^{\pm{\ify}})^*&:=&{\rm Hom}(\QQ^{\pm{\ify}},\QQ).
\end{eqnarray*}
We will write a linear form $\vp \in (\QQ^{+{\ify}})^*$ as
$\vp(\vec{x})=\sum_{k \geq 1} \vp_k x_k$ ($\vp_j\in \QQ$). 
Similarly, we write 
$\vp \in (\QQ^{-{\ify}})^*$ as
$\vp(\vec{x})=\sum_{k \leq -1} \vp_k x_k$ ($\vp_j\in \QQ$).

For the sequence 
$\io^+=(i_k)_{k \geq 1}$ (resp. $\io^-=(i_k)_{k \leq -1}$) and 
$k \geq 1$ (resp. $k \leq -1$), we set 
\[
\kp:={\rm min}\{l:l>k>0\,\, ({\rm resp.}\,0>l>k)\,\,
{\rm  and }\,\,i_k=i_l\},
\]
if it exists, and
$$
\km:={\rm max}\{l:0<l<k\,\,{\rm (resp.\,l<k<0)}\,\,
{\rm  and }\,\,i_k=i_l\},
$$
if it exists, otherwise $\kp=\km=0$.

We define a linear form $\beta_k$ $(k\geq0)$ on $\QQ^{+{\ify}}$ by
\begin{equation}
\beta_{k}(\vec x) := \begin{cases}
x_k+\sum_{k<j<\kp}\lan h_{i_k},\al_{i_j}\ran x_j+x_{\kp}
&
(k \geq 1),\\
0& (k = 0).
\end{cases}
\end{equation}
We also define a linear form $\beta_k$ $(k\leq 0)$ 
on $\QQ^{-{\ify}}$ by 
\begin{equation}
\beta_k(\vec y):= \begin{cases}
y_{\km}+\sum_{\km<j<k}\lan h_{i_k},\al_{i_j}\ran y_j+y_k
& 
(k \leq -1),\\
0 &(k=0).
\end{cases}
\end{equation}
By using these linear forms, let us
 define a piecewise-linear operator 
$S_k=S_{k,\io}$ on $(\QQ^{\pm{\ify}})^*$
as follows:
\begin{equation}
S_k(\vp):=\begin{cases}
\vp-\vp_k\beta_k & if\,\, \vp_k>0,\\
 \vp-\vp_k\beta_{\km} & if\,\, \vp_k\leq 0,\\
\end{cases}
\label{Sk}
\end{equation}
for $\vp(\vec x)=\sum \vp_k x_k\in (\QQ^{\pm{\ify}})^*$.
Here we set
\beqnn
\Xi_{\io^\pm} &:= &
\{S_{\pm j_l}\cd S_{\pm j_2}S_{\pm j_1}(\pm x_{j_0})\,|\,
l\geq0,j_0,j_1,\cd,j_l\geq1\},\\
\Sigma_{\io^\pm} & := &
\{\vec x\in \ZZ^{\pm{\ify}}\subset \QQ^{\pm{\ify}}
\,|\,\vp(\vec x)\geq0\,\,{\rm for}\,\,
{\rm any}\,\,\vp\in \Xi_{\io^\pm}\},
\eeqnn

We impose on $\io^+$ and $\io^-$ the following assumptions (P),(N):
$$
{\hbox{(P)   for $\io^+$, if a positive $k$ satisfies 
$\km=0$ then $\vp_k\geq0$ for any 
$\vp(\vec x)=\sum_k\vp_kx_k\in \Xi_{\io^+}$}},
$$
$$
{\hbox{(N)  for $\io^-$, if a negative $k$ satisfies 
$\kp=0$ then $\vp_k\leq0$ for any 
$\vp(\vec x)=\sum_k\vp_kx_k\in \Xi_{\io^-}$}}.
$$

\begin{thm1}[\cite{NZ}]
\label{nz}
Let $\io^{\pm}$ be the indices of sequences which are satisfied 
$(\ref{seq-con})$ and the assumptions (P),(N). 
Suppose $\Psi_{\io^+}:B(\ify)\hookrightarrow \ZZ^{\ify}_{\io^+}$ and
$\Psi_{\io^-}:B(-{\ify})\hookrightarrow \ZZ^{-{\ify}}_{\io^-}$ 
are the 
Kashiwara embeddings. Then, we have
${\rm Im}(\Psi_{\io^+})(\cong B(\ify))=\Sigma_{\io^+}$,\,
${\rm Im}(\Psi_{\io^-})(\cong B(-{\ify}))=\Sigma_{\io^-}$.
\end{thm1}

%%%%%%%%%% section 3  %%%%%%%%%%%%%%%
\section{Modified Quantum Algebras and its Crystal Base}
\setcounter{equation}{0}
\renewcommand{\thesection}{\arabic{section}}
\renewcommand{\theequation}{\thesection.\arabic{equation}}

%%%%%%%%%%% 3.1 %%%%%%%%%%%%%
%\subsection{Modified Quantum Algebra}

We define the left $\uq$-module $\uq a_{\lm}$ 
(\cite{K3}) by the relation: 
$q^h a_{\lm} = q^{\lan h,\lm\ran}a_{\lm}$.
Then $\util = \oplus_{\lm \in P}\uq a_{\lm}$ has 
an algebra structure by 
\begin{enumerate}
\item
$a_{\lm}P = Pa_{\lm-\xi}$ \,\, for \,\,$\xi \in Q$
and \,\,$P \in \uq_{\xi}$\\
$(\,\uq_{\xi}
:=\{P \in \uq;\,q^hPq^{-h}=q^{\lan h,\xi\ran}P\,
{\hbox{ for any }} h \in P^*\}\,)$.
\item
$a_{\lm}a_{\mu}=\del_{\lm,\mu}a_{\lm}$,
\end{enumerate}
and we call this algebra {\it modified quantum algebra}.

Let $M$ be a $\uq$-module with the weight space decomposition 
$M=\oplus_{\lm\in P}M_\lm$. Then $a_\lm$ is a projection 
$a_\lm:M\lar M_\lm$.

%%%%%%%%%%% 3.2 %%%%%%%%%%%%%
%\subsection{Crystal Base of $\uq a_{\lm}$}
In \cite{L}, it is revealed that 
modified quantum algebra $\util$
has a crystal structure and in \cite{K3} its
 crystal base is described as follows:

\begin{thm1}[\cite{K3}]
\begin{eqnarray*}
&&B(\uq a_{\lm})=B(\ify)\ot T_{\lm}\ot B(-{\ify}),\\
&&B(\util)=\bigoplus_{\lm\in P}B(\ify)\ot T_{\lm}\ot B(-{\ify}).
\end{eqnarray*}
\end{thm1}

%%%%%%%%%% section 4  %%%%%%%%%%%%%%%
\section{Polyhedral Realization of 
$B(\uq a_\lm)$}
\setcounter{equation}{0}
\renewcommand{\thesection}{\arabic{section}}
\renewcommand{\theequation}{\thesection.\arabic{equation}}

%%%%%%%%%%% 4.1 %%%%%%%%%%%%%
\subsection{Crystal structure of $\ZZ^{\ify}[\lm]$}

Let  $\ZZ^{+{\ify}}_{{\io}^+}$ and
 $\ZZ^{-{\ify}}_{{\io}^-}$ be as in Sect.2.2, 
We take the indices sequence $\io:=(\io^+,t_{\lm},\io^-)=
(\cd,i_2,i_1,t_{\lm},i_{-1},i_{-2},\cd)$ and weight $\lm\in P$.
We set
$\ZZ^{\ify}_{\io}[\lm]:=
\ZZ^{\ify}_{{\io}^+}\ot T_{\lm}\ot \ZZ^{-{\ify}}_{{\io}^-}.$
The crystal structure on 
$\ZZ^{\ify}_\io[\lm]$ associated with $\io$ and $\lm$ 
is defined as follows:
We identify $\ZZ^{+{\ify}} \ot T_{\lm} \ot \ZZ^{-{\ify}}$ with
$\ZZ^\ify$. Therefore, $\ZZ^{\ify}[\lm]$ is regarded as
 a sublattice of $\QQ^{\ify}=\QQ\ot_{\ZZ}\ZZ^\ify$. 
Thus, we can deonte $\vec x \in \ZZ^{\ify}_\io[\lm]$ by
$\vec x=(\cd,x_2,x_1,t_{\lm},x_{-1},x_{-2},\cd)$.
For $\vec x \in \QQ^{\ify}$, we define a linear function 
$\sigma_k(\vec x)\,\,(k \in \ZZ)$ by:
\beqn
\sigma_k(\vec x):= 
  \begin{cases}
    x_k+\sum_{j>k}\lan h_{i_k},\al_{i_j}\ran x_j
   &(k\geq1),\\
    -\lan h_{i_k},\lm\ran+x_k+\sum_{j>k}\lan h_{i_k},\al_{i_j}\ran x_j
   &(k\leq -1),\\
    -\ify
   &(k=0).
\end{cases}
\label{sigma}
\eeqn

Since $x_j=0$ for $j\gg0$,
 $\sigma_k$ is
well-defined.
Let $\sigma^{(i)} (\vec x)
 := {\rm max}_{k: i_k = i}\sigma_k (\vec x)$ and
\begin{equation}
M^{(i)} = M^{(i)} (\vec x) :=
\{k: i_k = i, \sigma_k (\vec x) = \sigma^{(i)}(\vec x)\}.
\label{m(i)}
\end{equation}

Note that 
$\sigma^{(i)} (\vec x)\geq 0$, and that 
$M^{(i)} = M^{(i)} (\vec x)$ is a finite set if and only if 
$\sigma^{(i)} (\vec x) > 0$. 
Now, we define the map 
$\eit: \ZZ^{\ify}[\lm] \sqcup\{0\}\lar \ZZ^{\ify}[\lm] \sqcup\{0\}$
, 
$\fit: \ZZ^{\ify}[\lm] \sqcup\{0\}\lar \ZZ^{\ify}[\lm] \sqcup\{0\}$ , by
 $\eit(0)=\fit(0)=0$ and 
\begin{equation}
(\fit(\vec x))_k  = x_k + \delta_{k,{\rm min}\,M^{(i)}}
\,\,{\rm if }\,\,M^{(i)} \,\,{\rm exists};
\,\,{\rm otherwise}\,\,\fit(\vec x)=0,
\label{action-f}
\end{equation}

\begin{equation}
(\eit(\vec x))_k  = x_k - \delta_{k,{\rm max}\,M^{(i)}} \,\, {\rm if}\,\,
M^{(i)} \,\,{\rm exists}; \,\,
 {\rm otherwise} \,\, \eit(\vec x)=0.
\label{action-e}
\end{equation}
where $\del_{i,j}$ is Kronecker's delta.
We also define the weight function and the function
$\vep_i$ and $\vp_i$ on $\ZZ^{\ify}[\lm]$ as follows:
\begin{equation}
\begin{array}{l}
wt(\vec x) :=\lm -\sum_{j=-{\ify}}^{\ify} x_j \al_{i_j}, \,\,
\vep_i (\vec x) := \sigma^{(i)} (\vec x),\\
\vp_i (\vec x) := \lan h_i, wt(\vec x) \ran + \vep_i(\vec x).
\end{array}
\label{wt-vep-vp}
\end{equation}

We denote this crystal by $\ZZ^{\ify}_{\io}[\lm]$.

Since there exist the embedings of crystals:
$B(\pm)\hookrightarrow \ZZ_{\io^\pm}^{\pm\ify}$, we obtain
\begin{thm1}
\begin{eqnarray*}
\Psi_{\io}^{(\lm)}:B(\ify)\ot T_{\lm}\ot B(-{\ify})
&\hookrightarrow&
\ZZ^{+{\ify}}_{{\io}^+}\ot T_{\lm}\ot \ZZ^{-{\ify}}_{{\io}^-}(=\ZZ^{\ify}_{\io}[\lm])\\
u_{\ify}\ot t_{\lm}\ot u_{-{\ify}} &\mapsto& (\cd,0,0,t_{\lm},0,0,\cd)
\end{eqnarray*}
is the unique stirict embedding which is associated with 
$\io:=(\cd,i_2,i_1,t_{\lm},i_{-1},i_{-2},\cd)$.
\label{Psi-lm}
 \end{thm1}

%%%%%%%%%%%% 4.2 %%%%%%%%%%%%%%%%
\subsection{The image of $\Psi^{(\lm)}_{\io}$}

Fix a sequence of indices $\io$ as above.
We define a linear function $\bar\beta_k(\vec x)$ as follows:
\beqn
\bar\beta_k (\vec x) & = & \sigma_k (\vec x) - \sigma_{\kp} (\vec x)
\label{beta}
\eeqn
where $\sigma_k$ is defined by (\ref{sigma}).
Since $\lan h_{i},\al_{i}\ran = 2$ for any $i \in I$, we have 
\beqnn
\bar\beta_k(\vec{ x}) & = &
  \begin{cases}
    x_k+\sum_{k<j<\kp}
    \lan h_{i_k},\al_{i_j}\ran x_j+x_{\kp}
    &(k\geq1\,\,or\,\,\kp \leq -1),\\
    -\lan h_{i_k},\lm\ran + 
    x_k+\sum_{k<j<\kp}
    \lan h_{i_k},\al_{i_j}\ran x_j+x_{\kp}
    &(k\leq-1\,\,and\,\,\kp > 0).
\end{cases}
\eeqnn

Using this notation, we define an operator $\bar{S_k}=\bar{S}_{k,\io}$
for a linear function $\vp(\vec x)=c+\sum_{-{\ify}}^{\ify}\vp_kx_k$
$(c,\vp_k\in\QQ)$ as follows:
\begin{equation}
 \bar{S_k}\,(\vp) :=\begin{cases}
 \vp - \vp_k \bar\beta_k & if\,\, \vp_k > 0,\\
 \vp - \vp_k \bar\beta_{k^{(-)}} & if\,\, \vp_k \leq 0.
\end{cases}
\label{S_k}
\end{equation}

An easy check shows $(\bar{S_k})^2=\bar{S_k}$.
For a sequence $\io$ and an integral weight $\lm$, we denote by
$\Xi_{\io}[\lm]$ the subset of linear forms which
 are obtained from 
the coordinate forms $x_j$, $x_{-j}$ ($j\geq1$) by applying 
transformations $\bar{S}_k$. In other words, we set
\begin{equation}
 \begin{array}{ll}
  \Xi^+_{\io}[\lm]&:=\{\bar{S}_{j_l}\cd \bar{S}_{j_1}(x_{j_0})\,
  :\,l\geq0,\,j_0,\cd,j_l>0\}
  \\
\Xi^-_{\io}[\lm]&:=
\{\bar{S}_{{-j}_k}\cd \bar{S}_{{-j}_1}(-x_{{-j}_0})\,
  :\,k\geq0,\,j_0,\cd,j_k>0\},\\
\Xi_\io[\lm]&:=\Xi^+_{\io}[\lm]\cup \Xi^-_{\io}[\lm].
 \end{array}
\label{Xi}
\end{equation}

Now we set 
\begin{equation}
 \Sigma_{\io}[\lm]
 :=\{\vec x\in \ZZ^{\ify}_{\io}[\lm](\subset \QQ^{\ify})\,:\,
\vp(\vec x)\geq 0\,\,{\rm for \,\,any }\,\,\vp\in \Xi_{\io}[\lm]\}.
 \label{Sigma}
\end{equation}
By Theorem \ref{nz}, we have
\begin{thm1}
 \label{main}
Suppose that $\io^{\pm}$ satisfys the assumption (P),(N),and 
(\ref {seq-con}).Let 
$\Psi^{(\lm)}_{\io}:
B(\ify)\ot T_{\lm}\ot B(-{\ify})\hookrightarrow 
\ZZ^{\ify}_{\io}[\lm]$
be the embedding of (\ref{Psi-lm}). Then 
${\rm Im}(\plm)(\cong B(\ify)\ot T_{\lm}\ot B(-{\ify}))$ 
is equal to 
$\Sigma_{\io}[\lm]$.
\end{thm1}

{\sl Remark.\,\,} 
Under the assumptions (P) and (N),  
both $\Xi^+_{\io}[\lm]$ and $\Xi^-_{\io}[\lm]$ are closed 
by the actions of $\bar S_k$'s, since 
$S_k$ $(k<0)$ (resp. $S_k$ $(k>0)$) acts identically on 
$\Xi^+_{\io}[\lm]$ (resp. $\Xi^-_{\io}[\lm]$).

%%%%%%%%%%%%%%%%%%%%%%%
The following lemma will be used lalter.

\begin{lem1}
\label{mlem}
Let $\Xi$ be a set of linear functions on $\QQ^{\ify}$.
Suppose that $\Xi$ 
is closed by actions of all $\bar{S}_k$'s, then 
the set 
\[
\Sigma=\{\vec x\in \ZZ^{\ify}_{\io}[\lm]|
\vp(\vec x)\geq0{\hbox{ for any }}\vp\in \Xi_\io\}
\] 
is a sub-crystal of $\ZZ^\ify_\io[\lm]$.
 \end{lem1}

{\sl Proof.}\,\,
It suffices to show that $\Sigma$ 
is closed under the actions of $\fit$ and $\eit$.
For 
$\vec x \in \Sigma$, suppose  
$\fit\vec x = (\cd,x_{k}+1,\cd)$. For any $\vp =
 c+\sum\vp_jx_j\in \Xi\,\,(c,\vp_j \in \QQ)$, we need to
show that
\begin{equation}
\vp(\fit\vec x) \geq 0.
\end{equation}

Since $\vp(\fit(\vec x)) = \vp(\vec x) + \vp_{k} \geq \vp_{k}$,
it is enough to consider the case when $\vp_{k} <0$. 
By definition of $\fit\vec x$, we have 
$\sigma_{\km} < \sigma_{k}$. This shows that
\begin{eqnarray*}
\sigma_{\km} < \sigma_{k} 
&\Longleftrightarrow& \bar\beta_{\km}
=\sigma_{\km} - \sigma_{k} < 0 \\
&\Longrightarrow& \bar\beta_{\km} \leq -1. 
\end{eqnarray*}
Therefore, it follows that
\begin{eqnarray*}
\vp(\fit\vec x) &=& \vp(\vec x) + \vp_{k} \\
                &\geq& \vp(\vec x) - \vp_{k}\bar\beta_{\km} \\
                &=& (\bar{S}_{k}\vp)(\vec x) \,\, \geq 0.
\end{eqnarray*}

Suppose that $\eit\vec x = (\cd,x_{k}-1,\cd)$.
We need to show that
\begin{equation}
\vp(\eit\vec x) \geq 0.
\end{equation}

Since $\vp(\eit(\vec x)) = \vp(\vec x) - \vp_{k} \geq -\vp_{k}$,
it is enough to consider the case when $\vp_{k} >0$.
By definition of $\eit\vec x$, we have 
$\sigma_{k} > \sigma_{\kp}$.
This shows that 
\begin{eqnarray*}
 \sigma_{k} > \sigma_{\kp} 
&\Longleftrightarrow
& \bar\beta_{k}=\sigma_{k} - \sigma_{\kp} > 0 \\
&\Longrightarrow& \bar\beta_{\km} \geq 1. 
\end{eqnarray*}
Therefore, it follows that
\begin{eqnarray*}
\vp(\eit\vec x) &=& \vp(\vec x) - \vp_{k} \\
                &\geq& \vp(\vec x) - \vp_{k}\bar\beta_{k} \\
                &=& (\bar{S}_{k}\vp)(\vec x) \,\,\geq 0.
\end{eqnarray*}
\qed

%%%%%%%%% Section 5  %%%%%%%%%%
\section{Polyhedral Realization of $B(\uq a_{\lm})$ of Type $A_n$}
\setcounter{equation}{0}
\renewcommand{\theequation}{\thesection.\arabic{equation}}
%%%%%%%%%%%%%%%%%%%%
In this section, we shall describe the crystal structure 
of the component including $t_\lm$
in $B(\ify)\ot T_{\lm}\ot B(-{\ify})$
for the case of type $A_n$.

It will be convenient for us to change the indexing set for
$\ZZ^{\ify}$ from $\ZZ_{\geq 1}$ to 
$\ZZ_{\geq 1} \times [1,n]$. We will do this with the help of the 
bijection $\ZZ_{\geq1} \times [1,n] \to \ZZ_{\geq 1}$ given by  
($(j;i) \mapsto (j-1)n + i$). Thus, we will write an element
 $\vec x \in \ZZ^{+{\ify}}$ 
as doubly-indexed family $(x_{j;i})_{j \geq 1, i \in [1,n]}$ of
nonnegative integers. Simillarly,
using that
$\ZZ_{\geq 1} \times [1,n] \to \ZZ_{\leq -1}$ 
($(j;i) \mapsto -jn + i - 1$) is bijective, we will write an element
 $\vec x \in \ZZ^{-{\ify}}$ 
as doubly-indexed family $(x_{-j;i})_{j \geq 1, i \in [1,n]}$ of 
nonpositive integers.
Therefore, we can write that $\vec x \in \ZZ^{\ify}$ as   
$(\cd,x_{1;2},x_{1;1},t_{\lm},x_{-1;n},x_{-1;n-1},\cd)$.
We will adopt the convention that $x_{j;0} = x_{j;n+1} = 0$
unless $i \in [1,n]$

To state the main theorem, we prepare several things.
For $x \in \mathbb R$, set $(x)_+ := max (0,x)$.
Let $\lm=\sum_{1\leq i\leq n}\lm_i\Lm_i$ be an integral weight
satisfying 
 $\lm_1,\cd,\lm_{i_0}>0$ and 
$\lm_{{i_0}+1},\cd,\lm_{n}\leq0$ for some $i_0$ and 
for $(j;i)\in \ZZ_{\geq1}\times [1,n]$ set 
\[
C_{-j;i}=
\left\{\begin{array}{ll}
(-{\lm_{-j+i+1}}+(-{\lm_{-j+i+2}}+( \cd + (-{\lm_{-j+n+1}})_+)_+ \cd
)_+)_+&{\hbox{if }}1\leq j\leq i\leq n,\\
0&{\hbox{otherwise}}.
\end{array}\right.
\]
We will use the following lemma frequently:
\begin{lem1} 
\label{plem}
For real numbers  $r_1,\cd, r_n,$  we have,
\begin{equation*}
r_1+(r_2+(r_3+\cd+(r_{n-1}+(r_n)_+)_+)_+)_+
={\rm max}(r_1,r_1+r_2,\cd,r_1+r_2+\cd+r_n)
\label{r1-n}
\end{equation*}
\end{lem1}
{\sl Proof.}

We can easily show from the fact :
$r_1+(r_2)_+={\rm max}(r_1,r_1+r_2)$ and iterating this.
\qed
\q\\

By the above lemma, we can write
$$C_{-j;i}={\rm max}(0,-\lm_{-j+i+1},-\lm_{-j+i+1}-\lm_{-j+i+2},\cd,
-\lm_{-j+i+1}-\lm_{-j+i+2}-\cd -\lm_{-j+n+1}). $$
\begin{thm1} 
 \label{A_n}
Let $\io = (\cd,2,1,n,\cd,2,1,t_{\lm},n,n-1,\cd,1,n,n-1,\cd)$ be an
infinite sequence and 
$\lm$ and $C_{-j;i}$ be as above. 
We define 
\begin{eqnarray*}
\Xi^{'}_{\io}[\lm]&:=&
\{\bar{S}_{-{j_k}}\cd \bar{S}_{-{j_1}}({x_{{-j};i}}+C_{-j;i})\,
:\,k\geq0,\,i\in I,j\geq1,\,\,j_1,\cd,j_k\geq1\},\\
\Sigma'_\io[\lm]&:=&
\{\vec x\in \ZZ^{\ify}_{\io}[\lm](\subset \QQ^{\ify})\,:\,
    \vp(\vec x)\geq 0\,\,{\rm for \,\,any }\,\,\vp\in
 \Xi^{'}_{\io}[\lm]\}
\end{eqnarray*}
and denote the connected component of ${\rm Im} \,
 (\Psi^{(\lm)}_{\io })$ containing $\vec 0 := 
(\cd,0,0,t_{\lm},0,0,\cd)$ by $B_0(\lm)$.
Then we have 
\begin{enumerate}
\item
$B_0(\lm) = \Sigma_{\io}[\lm]\cap\Sigma'_{\io}[\lm]$.
\item
Let $v_{\lm}$ be the unique highest weight vector in
 $B_0(\lm)$. Then we have
\[
v_{\lm}
= (\cd,0,0,t_{\lm},-C_{-1;n},-C_{ -1;n-1},\cd,-C_{-j;i},\cd),
\]
%$($in the case that $1 \leq i < j$, by $(\cite{NZ})$, we know
%$x_{-j;i} =0)$.
\end{enumerate}
\label{thm1}
\end{thm1}

{\sl Proof.}\,\,
Since $\Xi^{'}_{\io}[\lm]$ 
is closed under the actions of $\bar{S}_{k}$'s,
by Lemma \ref{mlem}
$\Sigma_\io[\lm]\cap \Sigma'_\io[\lm]$ 
has a crystal structure unless it is empty. 
We will show that $\Sigma_\io[\lm]\cap \Sigma'_\io[\lm]$
 contains $\vec 0$,
which implies that $\Sigma_\io[\lm]\cap \Sigma'_\io[\lm]$
 is non-empty,
and has the unique highest weight vector. 
First, we will show that 
\begin{equation}
   x_0 =(\cd,0,0,t_{\lm},-C_{-1;n},
    -C_{ -1;n-1},\cd,-C_{-j;i},\cd)\,\,
\label{x0}
%   (1\le j\le i \le n)
\end{equation}
is a highest weight vector in $\ZZ^\ify_\io[\lm]$.
We set
$y_{-j;i}:=-{\lm_{-j+i+1}}+(-{\lm_{-j+i+2}}+(\cd + 
    (-{\lm_{-j+n+1}})_+)_+)_+$. Thus, we have 
$C_{-j;i} = -(y_{-j;i})_+$.
Due to (\ref{action-e}), it suffices to show 
$$
\sigma_{-j:i}(x_0)\leq0 {\hbox{ for }}j\geq1, \,\,i\in I.
$$
(In the case $j<0$, trivially $\sigma_{-j:i}(x_0)=0$.)
We consider the following four cases:

(I) $j=1$. \,\,(II) $i=n$. \,\,(III) $1\leq i < j\leq n$.\,\,
(IV) $1<j\leq i <n$.

\vskip5pt
(I) The case $j=1$.

\noindent
We will show that $\sigma_{-1,n}(x_0),\sigma_{-1;n-1}(x_0),\cd
,\sigma_{-1;1}(x_0) \leq 0 $.
Note the following simple fact:

\begin{equation}
-(-a)_+ -a \leq 0,
\q {\hbox{for any }}a \in \mathbb R.
\label{zeroika}
\end{equation}

We can write $y_{-1;i+1}:= -\lm_{i+1}+(-\lm_{i+2}+(\cd+(-\lm_n)_+)_+\cd)_+$.
By the definition of $\sigma_{-j:i}$, we have 
$\sigma_{-1;i}(x_0) = -(-\lm_{i} + (y_{-1;i+1})_+)_+
 + (y_{-1;i+1})_+ - \lm_i$.
By $(\ref{zeroika})$, we obtain $\sigma_{-1;i}(x_0)\leq 0$.
This shows 
$\sigma_{-1;n}(x_0),\sigma_{-1;n-1}(x_0),\cd,\sigma_{-1;1}(x_0) \leq 0$.

\vskip7pt
(II) The case $i=n$.

\noindent
We shall show  
$\sigma_{-j;n}(x_0) \leq 0 \,\,(1 \leq \forall j \leq n)$ 
by the induction on $j$. If $j=1$, it is true by (I). 
Suppose $j>1$. We can write $\sigma_{-j;n}(x_0) =
-(-\lm_{-j+n+1})_+ + (-\lm_{-j+n+1} + (-\lm_{-j+n+2})_+)_+
-(-\lm_{-j+n+2})_+  + \sigma_{-j+1;n}(x_0)$.
Now, set $A:=\sigma_{-j;n}(x_0) - \sigma_{-j+1;n}(x_0)$.
Since $\sigma_{-j+1;n}(x_0) \leq 0$ by the induction hypothesis, 
it is sufficient to show $A\leq0$.
%If $-\lm_{-j+n+1} \geq 0$, then
%$A=\lm_{-j+n+1} + (-\lm_{-j+n+1} + (-\lm_{-j+n+2})_+)_+
%-(-\lm_{-j+n+2})_+$. By $\ref{zeroika}$

If $\lm_{-j+n+1}\, , \,\lm_{-j+n+2} \geq 0$ 
$\lm_{-j+n+1}\, , \,\lm_{-j+n+2} \leq 0$, then obviously 
$A = 0$. If $\lm_{-j+n+1} \geq 0$ and $\lm_{-j+n+2} \leq 0 $,
we can write $A = (-\lm_{-j+n+1} - \lm_{-j+n+2})_+ + \lm_{-j+n+2}$.
In this case, if $-\lm_{-j+n+1} - \lm_{-j+n+2} \leq 0$, then 
$A = \lm_{-j+n+2} \leq 0$.
If $-\lm_{-j+n+1} - \lm_{-j+n+2} \geq 0$, then $A = -\lm_{-j+n+1} \leq 0$.
If $\lm_{-j+n+1} \leq0$ and $\lm_{-j+n+2}\geq0$, obviously $A=0$.

\vskip7pt
(III) The case $1\leq i < j \leq n$.

\noindent
By the definition, $C_{-j;i}=0$ for $i<j$. We can write
$$
\sigma_{-j;i}(x_0)=C_{-i-1:i+1}-C_{-i:i}+\sigma_{-i;i}(x_0).$$
By Lemma $\ref{r1-n}$, we have
\begin{align*}
C_{-i:i}&={\rm max}(0,-\lm_{1},-\lm_{1}-\lm_{2},\cd,
-\lm_{1}-\lm_{2}-\cd -\lm_{-i+n+1}),\\
C_{-i-1:i+1}&={\rm max}(0,-\lm_{1},-\lm_{1}-\lm_{2},\cd,
-\lm_{1}-\lm_{2}-\cd -\lm_{-i+n}).
\end{align*}
Therefore, we obtain $C_{-i:i}\geq C_{-i-1:i+1}$.
This shows $\sigma_{-j;i}(x_0)\leq 0$.

\vskip7pt
(IV) The case $1<j\leq i <n$.

\nd
We will show $\sigma_{-j;i}(x_0)\leq0$ by the induction on
$1 < j\leq i < n$. 
We can write:
\begin{eqnarray}
\sigma_{-j;i}(x_0)&=&-C_{-j;i}+C_{-j;i+1}+C_{-j+1;i-1}-C_{-j+1;i}+
\sigma_{-j+1;i}(x_0)\nn
\end{eqnarray}
where
\begin{eqnarray*}
C_{-j;i}&=&(-{\lm_{-j+i+1}}+(-{\lm_{-j+i+2}}+ \cd + 
    (-{\lm_{-j+n+1}})_+\cd)_+)_+,\\
C_{-j;i+1}&=&(-{\lm_{-j+i+2}}+(-{\lm_{-j+i+2}}+ \cd + 
    (-{\lm_{-j+n+1}})_+\cd)_+)_+,\\
C_{-j+1;i-1}&=&(-{\lm_{-j+i+1}}+(-{\lm_{-j+i+2}}+ \cd +
(-{\lm_{-j+n+1}}+(-{\lm_{-j+n+2}})_+)_+\cd)_+)_+,\\
C_{-j+1;i}&=&(-{\lm_{-j+i+2}}+(-{\lm_{-j+i+3}}+ \cd + 
    (-{\lm_{-j+n+1}}+(-\lm_{-j+n+2})_+)_+\cd)_+)_+.
\end{eqnarray*}
Set 
$B:=\sigma_{-j;i}(x_0) - \sigma_{-j+1;i}(x_0)$. 
It is sufficient to show $B\leq0$, since
 $\sigma_{-j+1;i}(x_0) \leq 0$ by the induction hypothesis.
By Lemma $\ref{r1-n}$, we have
\begin{eqnarray*}
C_{-j;i+1}&=&{\rm max}(0,-\lm_{-j+i+2},-\lm_{-j+i+2}-\lm_{-j+i+3},\cd,
-\lm_{-j+i+2}-\lm_{-j+i+3}-\cd -\lm_{-j+n+1}),\\
C_{-j+1;i}&=&{\rm max}(0,-\lm_{-j+i+2},-\lm_{-j+i+2}-\lm_{-j+i+3},\cd,
-\lm_{-j+i+2}-\lm_{-j+i+3}-\cd -\lm_{-j+n+2}),
\end{eqnarray*}
and then 
\begin{eqnarray}
 C_{-j;i+1}\leq C_{-j+1;i}.
\label{cc}
\end{eqnarray}
Since
\begin{eqnarray*}
B&=&-(-{\lm_{-j+i+1}}+C_{-j;i+1})_+ +C_{-j;i+1} +
(-{\lm_{-j+i+1}}+C_{-j+1;i})_+ -C_{-j+1;i},
\end{eqnarray*}
if $C_{-j;i+1}=C_{-j+1;i}$, then $B=0$. We consider the case 
 $C_{-j;i+1} \ne C_{-j+1;i}$. In this case, we have 
$C_{-j;i+1} < C_{-j+1;i}$ by $(\ref{cc})$.
We can write
\begin{eqnarray}
C_{-j;i+1}&=&0 \,\text{ or }-\lm_{-j+i+2}-\lm_{-j+i+3}-\cd -\lm_{-j+k+1}>0
\,\text{ for some }k \,( i \leq k \leq n),\nn \\
C_{-j+1;i}&=&-\lm_{-j+i+2}-\lm_{-j+i+3}-\cd -\lm_{-j+n+2}>0,\nn \\
C_{-j;i+1}-C_{-j+1;i}&=&\lm_{-j+k+2}+\lm_{-j+k+3}+\cd +\lm_{-j+n+2}<0.
\label{c_ji}
\end{eqnarray}
If $C_{-j;i+1}=0$, then
\begin{eqnarray*}
B&=&-(-{\lm_{-j+i+1}})_++(-{\lm_{-j+i+1}}+C_{-j+1;i})_+ -C_{-j+1;i}\\
&=&-(-{\lm_{-j+i+1}})_++(-\lm_{-j+i+1}-\lm_{-j+i+2}-\cd-\lm_{-j+n+2})_+\\
&&+\lm_{-j+i+2}+\lm_{-j+i+3}+\cd +\lm_{-j+n+2}.
\end{eqnarray*}
In this case, if 
$-\lm_{-j+i+1}-\lm_{-j+i+2}-\cd-\lm_{-j+n+2} \leq 0$, then 
$-\lm_{-j+i+1} \leq 0$ by $(\ref{c_ji})$. Therefore, we have
$B<0$. If
$-\lm_{-j+i+1}-\lm_{-j+i+2}-\cd-\lm_{-j+n+2} > 0$, then
$B=-(-{\lm_{-j+i+1}})_+ -{\lm_{-j+i+1}} \leq 0$ by $(\ref{zeroika})$.
Therefore, in the case $C_{-j;i+1}=0$, We get $B\leq 0$.
  
We consider the case $C_{-j;i+1}\ne 0$.
In this case, we have
\begin{eqnarray*}
B&=&-(-{\lm_{-j+i+1}}-\lm_{-j+i+2}-\cd -\lm_{-j+k+1})_+\\
&&+(-\lm_{-j+i+1}-\lm_{-j+i+2}-\cd-\lm_{-j+k+1}-
\lm_{-j+k+2}-\cd -\lm_{-j+n+2})_+\\
&&+\lm_{-j+k+2}+\lm_{-j+k+3}+\cd +\lm_{-j+n+2}.
\end{eqnarray*}
If $-{\lm_{-j+i+1}}-\lm_{-j+i+2}-\cd -\lm_{-j+k+1} \geq 0$, 
 by $(\ref{c_ji})$, we have
\begin{eqnarray*}
B&=&{\lm_{-j+i+1}}+\lm_{-j+i+2}+\cd +\lm_{-j+k+1}\\
&&-\lm_{-j+i+1}-\lm_{-j+i+2}-\cd-\lm_{-j+k+1}-
\lm_{-j+k+2}-\cd -\lm_{-j+n+2}\\
&&+\lm_{-j+k+2}+\lm_{-j+k+3}+\cd +\lm_{-j+n+2}\\
&=&0.
\end{eqnarray*}
On the other hand, if 
$-{\lm_{-j+i+1}}-\lm_{-j+i+2}-\cd -\lm_{-j+k+1} < 0$, by
$(\ref{c_ji})$, we have
\begin{eqnarray*}
B&=&(-\lm_{-j+i+1}-\lm_{-j+i+2}-\cd-\lm_{-j+k+1}-
\lm_{-j+k+2}-\cd -\lm_{-j+n+2})_+\\
&&+\lm_{-j+k+2}+\lm_{-j+k+3}+\cd +\lm_{-j+n+2}\\
&<&0.
\end{eqnarray*}
This shows $B \leq 0$.

Therefore, we conclude  
$\sigma_{-j;i}(x_0) = B + \sigma_{-j+1;i}(x_0)\leq 0$ for any $j,i$ 
$(1\leq i,j\leq n)$ by the induction hypothesis.
Now, we have shown that 
$x_0=(\cd,0,0,t_{\lm},-C_{-1;n},
   -C_{ -1;n-1},\cd,-C_{-j;i},\cd)$ 
is one of the highest weight vecters in $\ZZ^\ify_{\io}[\lm]$.
\vskip7pt

Next, we will show that 
$\Sigma_{\io}[\lm]
    \cap
\Sigma'_{\io}[\lm]$
contains $\vec 0$.

We know that $\vec 0 \in \Sigma_{\io}[\lm]$, since 
all the constant terms of all elements in $\Xi_{\io}[\lm]$ are $0$. 
Therefore, we will show  $\vec 0 \in\Sigma'_{\io}[\lm]$.

For the purpose, we shall evaluate the constant term of 
$\bar{S}_{-{j_k}}\cd \bar{S}_{-{j_1}}({x_{{-j};i}})$ 
by using the similar method to the one in (\cite{NZ}, Lemma 5.2)
and see that all of them are non-negative.

For $l\geq 0$, we define $\bar{S}^{(l)}_{-j;i}$ as follows: 
\begin{equation}
\bar{S}^{(l)}_{-j;i}:= \bar{S}_{-j;i+l-1}\bar{S}_{-j;i+l-2}
\cd\bar{S}_{-j;i+1}\bar{S}_{-j;i}.
\label{const1}
\end{equation}
Now, we define 
\[
M^{(i)}:=\{\mu=(\mu_1,\mu_2,\cd,\mu_i)\in \ZZ^i|
n-i+1\geq \mu_1\geq\mu_2\geq\cd\geq\mu_i\geq0\},
\]
whose element is called {\it $i$-admissible partition} (\cite{NZ}).
For $\mu\in M^{(i)}$ set
\begin{equation}
\vp^{(\mu)}_{-j;i}(x)
:= \bar{S}^{(\mu_{i})}_{-j+i-1;1}
\bar{S}^{(\mu_{i-1})}_{-j+i-2;2}\cd
\bar{S}^{(\mu_{j+1})}_{1;i-j}
\bar{S}^{(\mu_{j})}_{-1;i-j+1}
\bar{S}^{(\mu_{j-1})}_{-2;i-j+2}\cd\bar{S}^{(\mu_{2})}_{-j+1;i-1}
\bar{S}^{(\mu_{1})}_{-j;i}(x_{-j;i})
\label{const2}
\end{equation}
By the similar arguement in \cite{NZ} Lemma 5.2, we obtain
the explicit form of $\vp^{(\mu)}_{-j;i}(x)$ up to constant term:
\begin{lem1}
\label{vpij}
We have 
\begin{eqnarray}
&& \Xi'_\io[\lm]=
\{\vp^{(\mu)}_{-j;i}(x)+C_{-j;i}|
j\geq1,\,i\in I,\,\mu\in M^{(i)}\},\\
&&\vp^{(\mu)}_{-j;i}(x)-\vp^{(\mu)}_{-j;i}(0)=
\sum_{k=1}^{i}(x_{-j+k-\theta(j-k);i-k+1+\mu_{k}}
-x_{-j+k+1-\theta(j-k-1);i-k+\mu_{k}}),
\end{eqnarray}
where we consider $x_{j:i}=0$ unless $1\leq i\leq n$ and 
\[
\theta(x):=
\begin{cases}
1&\text{if $x\geq0$,}\\
0&\text{if $x<0$.}
\end{cases}
\]
\end{lem1}
Let us evaluate the constant term in $\vp^{(\mu)}_{-j;i}$.
For $\mu=(\mu_1,\cd,\mu_i)$, we set 
$\mu':=(\mu_1,\cd,\mu_{j-1},0,\cd,0)$ and 
$\mu'':=(\mu_1,\cd,\mu_{j-1},\mu_j,0,\cd,0)$
Then we know that 
$\vp^{(\mu')}_{-j;i}$ has no constant term by its definition. 
Thus, be Lemma \ref{vpij} we have
\begin{eqnarray*}
\vp^{(\mu)}_{-j;i} 
&=& (x_{-j;i+{\mu_1}} - x_{-j+1;i-1+{\mu_1}})+
(x_{-j+1;i-1+{\mu_{2}}} - x_{-j+2;i-2+{\mu_2}})
+ \cd \\
&& +(x_{-2;i-j+2+{\mu_{j-1}}} - x_{-1;i-j+1+{\mu_{j-1}}})
 +x_{-1;i-j+1},
\end{eqnarray*}
By direct calculations, 
we obtain 
%\begin{eqnarray*}
%\bar{S}^{(1)}_{-1;i-j+1}(\vp^{(\mu)}_{-j;i})
%&=&\bar{S}_{-1;i-j+1}(\vp^{(\mu)}_{-j;i})\\
%&=&\vp^{(\mu)}_{-j;i} -x_{-1;i-j+1} +x_{-1;i-j+2} + \lm_{i-j+1} 
%\end{eqnarray*}
\begin{eqnarray*}
\vp^{(\mu'')}_{-j;i}=
\bar{S}^{(\mu_j)}_{-1;i-j+1}(\vp^{(\mu')}_{-j;i})
&=& \bar{S}_{-1;i-j+\mu_j}
\cd\bar{S}_{-1;i-j+1}(\vp^{(\mu)}_{-j;i})\\
&=&\vp^{(\mu)}_{-j;i} -x_{-1;i-j+\mu_j} +x_{-1;i-j+\mu_j+1} 
+ \lm_{i-j+1} + \lm_{i-j+2}+\cd+\lm_{i-j+\mu_j}.
\end{eqnarray*}
Since $\bar{S}^{(\mu_{i})}_{-j+i-1;1}
\bar{S}^{(\mu_{i-1})}_{-j+i-2;2}\cd
\bar{S}^{(\mu_{j+1})}_{1;i-j}$ does not produce 
non-trivial constant term, we have 
\[
 \vp^{(\mu'')}_{-j;i}(0)=\vp^{(\mu)}_{-j;i}(0).
\]
Hence, we obtain the constant term of $\vp^{(\mu)}_{-j;i}(x)$, 
\[
 \vp^{(\mu)}_{-j;i}(0)=
\lm_{i-j+1}+\lm_{i-j+2}+\cd+\lm_{i-j+\mu_j}.
\]
By Lemma \ref{plem}, we have 
\[
\begin{array}{l}
C_{-j;i}=(-{\lm_{-j+i+1}}+(-{\lm_{-j+i+2}}
+ (\cd + (-{\lm_{-j+n+1}})_+)_+ )_+)_+\\
={\rm max}(0,-\lm_{-j+i+1},-\lm_{-j+i+1}-\lm_{-j+i+2},\cd,
-\lm_{-j+i+1}-\cd-\lm_{n-j+1}),
\end{array}
\]
and then 
\begin{equation}
 \begin{cases}
C_{-j;i}+\lm_{i-j+1} \geq 0,\\
 C_{-j;i}+\lm_{i-j+1}+\lm_{i-j+2} \geq 0,\\
\cd,\\
C_{-j;i}+\lm_{i-j+1}+\lm_{i-j+2}+\cd+\lm_{n-j+1} \geq 0,
 \label{csum}
 \end{cases}
\end{equation}
which implies  that $C_{-j;i}+\vp^{(\mu)}_{-j;i}(0)\geq0$.
Therefore, 
constant terms in all elements in $\Xi'_{\io}[\lm]$ 
are non-negative
and then 
we know that  
$\vec 0$ is conatained in  $\Sigma_{\io}[\lm]
    \cap \Sigma^{'}_{\io}[\lm]$.
Therefore, 
we obtain $B_0(\lm)\subset \Sigma_\io[\lm]\cap \Sigma'_\io[\lm]$.

Finally, we will show that $\Sigma_\io[\lm]\cap \Sigma'_\io[\lm]$ 
has the unique highest weight vector, which implies 
the opposite inclusion. We also show that 
the image of the highest weight vector 
by $\Psi_\io^{(\lm)}$ coincides with $v_{\lm}$
as in Theorem \ref{A_n}.

Let $x_0$ be the highest weight vector as in (\ref{x0}).
We set 
\[
 D_{-j;i}:=\sigma_{-j;i}(x_0) + C_{-j;i}.
\]

We shall show the following lemma on $D_{-j;i}$:
\begin{lem1}
\label{lem-d}
\begin{enumerate}
 \item 
If $D_{-j;i}\leq 0$, then we have $C_{-j;i}=0$.
\item
If $D_{-j;i}> 0$, then we have 
$$D_{-j;i}=
(-\lm_{-j+i+1}+(-\lm_{-j+i+2}
+(\cd+(-\lm_{-j+n+1})_+)_+)_+)_+ = C_{-j;i}.$$
\end{enumerate}

\end{lem1}

{\sl Proof.}

(I) The case $j=1$. \\
\noindent
By Lemma \ref{plem}, we get 
\begin{eqnarray*}
D_{-1;i}&=& (-\lm_{i+1}+(-\lm_{i+2}+(\cd +(-\lm_{n})_+)_+\cd )_+ -\lm_{i}\\
&=&{\rm max}(-\lm_i,-\lm_i-\lm_{i+1},\cd,-\lm_i-\lm_{i+1}-\cd-\lm_n),\\
C_{-1;i}&=&{\rm max}
(0,-\lm_i,-\lm_i-\lm_{i+1},\cd,-\lm_i-\lm_{i+1}-\cd-\lm_n)\\
&=&{\rm max}(0,D_{-1;i})
\end{eqnarray*}
If $D_{-1;i} \leq 0$, obviously $C_{-1;i}=0$.
If $D_{-1;i} > 0$, we obtain $D_{-1;i}=C_{-1;i}$
This shows that lemma is true.
\vskip7pt

(II) The case that $i=n$.

\noindent
We will show by the induction on $j$.
By the definition $D_{-j;i}$ and $C_{-j;i}$, we have
\begin{eqnarray}
D_{-j;n}&=&(-\lm_{-j+n+1}+(-\lm_{-j+n+2})_+)_+
-(-\lm_{-j+n+2})_++\sigma_{-j+1;n}(x_0),\\
C_{-j;n}&=&(-\lm_{-j+n+1})_+. \nn
\label{i=n0}
\end{eqnarray}

(i) The case $D_{-j+1;n}>0$.

\noindent
In this case, by the induction hypothesis, we have 
$\sigma_{-j+1;i}(x_0)=0$.
Then, 
\begin{eqnarray}
D_{-j;n}=(-\lm_{-j+n+1}+(-\lm_{-j+n+2})_+)_+-(-\lm_{-j+n+2})_+.
\label{i=n1}
\end{eqnarray}
\vskip5pt
(i-a) The case $D_{-j;n}\leq 0$.

\noindent
If $-\lm_{-j+n+2}\leq 0$, then $-\lm_{-j+n+1}\leq 0$.
Thus, $C_{-j;n}=0$. If $-\lm_{-j+n+2} > 0$, then 
$-\lm_{-j+n+1}\leq 0$ since $D_{-j;n}\leq 0$.
Therefore, we obtain $C_{-j;n}=0$. 

\vskip5pt
(i-b) The case $D_{-j;n} > 0$.

\noindent
Since the right hand-side of $(\ref{i=n1})$ is positive, we have
\begin{eqnarray}
{\rm max}(0, -\lm_{-j+n+1}, -\lm_{-j+n+1}-\lm_{-j+n+2})>
{\rm max}(0, -\lm_{-j+n+2}).
\label{i=n2}
\end{eqnarray}
This shows that the left hand-side of $(\ref{i=n2})$ is positive
and equal to $-\lm_{-j+n+1}$ or $-\lm_{-j+n+1}-\lm_{-j+n+2}$.
If L.H.S of $(\ref{i=n2})=-\lm_{-j+n+1}>0$, then
$-\lm_{-j+n+2} \leq 0$. Therefore, we obtain
$D_{-j;n}=-\lm_{-j+n+1}=C_{-j;n}$.
If L.H.S of $(\ref{i=n2})=-\lm_{-j+n+1}-\lm_{-j+n+2}>0$, then
$-\lm_{-j+n+2} > 0$. Then
$D_{-j;n}=-\lm_{-j+n+1}$ must be positive. Thus, we obtain 
$C_{-j;n}=-\lm_{-j+n+1}=D_{-j;n}$.
\vskip5pt
(ii) The case $D_{-j+1;n} \leq 0$.

\noindent
By the induction hypothesis, we have 
$C_{-j+1;n}=(-\lm_{-j+n+2})_+=0$ and then,
 by the definition of $\lm$, $\lm_{-j+n+1} \geq 0$.
Then, by $(\ref{i=n0})$, we have
\begin{eqnarray*}
D_{-j;n}&=&\sigma_{-j+1;n}(x_0) \leq 0.
\end{eqnarray*}
In this case, we have $C_{-j;n}=(-\lm_{-j+n+1})_+=0$.
\vskip5pt

(III) The case $1 \leq i < j \leq n$.

\noindent
By the definition, $C_{-j;i}=0$. Then, we have
$$D_{-j;i}=\sigma_{-j;i}(x_0)\leq 0.$$
In this case, the lemma is true.
\vskip5pt
(IV) The case $1 < j \leq i < n$.

\noindent
We have
\begin{eqnarray}
D_{-j;i}&=&C_{-j;i+1}+C_{-j+1;i-1}-C_{-j+1;i}+\sigma_{-j+1;i}(x_0).
\end{eqnarray}

\vskip5pt
(i) The case $D_{-j+1;i} > 0 $.

\noindent
By the induction hypothesis, we have
\begin{eqnarray*}
D_{-j+1;n}&=&(-\lm_{-j+i+2}+(-\lm_{-j+i+3}+\cd+(-\lm_{-j+n+2})_+\cd)_+)_+
=C_{-j+1;i},
\end{eqnarray*}
and then $\sigma_{-j+1;n}(x_0)=D_{-j+1;n}-C_{-j+1;n}=0$.
Then, we have
\begin{eqnarray*}
D_{-j;i}&=&(-\lm_{-j+i+2}+(-\lm_{-j+i+3}+\cd+(-\lm_{-j+n+1})_+)\cd_+)_+\\
&&+(-\lm_{-j+i+1}+(-\lm_{-j+i+2}+\cd+(-\lm_{-j+n+2})_+)\cd_+)_+\\
&&-(-\lm_{-j+i+2}+(-\lm_{-j+i+3}+\cd+(-\lm_{-j+n+2})_+)\cd_+)_+\\
&=&{\rm max}(0,-\lm_{-j+i+2},-\lm_{-j+i+2}-\lm_{-j+i+3},\cd,
-\lm_{-j+i+2}-\lm_{-j+i+3}-\cd -\lm_{-j+n+1})\\
&&+(-{\lm_{-j+i+1}}+{\rm max}(0,-\lm_{-j+i+2},-\lm_{-j+i+2}-\lm_{-j+i+3},\cd,
-\lm_{-j+i+2}-\lm_{-j+i+3}-\cd -\lm_{-j+n+2}))_+\\
&&-{\rm max}({\rm max}(0,-\lm_{-j+i+2},-\lm_{-j+i+2}-\lm_{-j+i+3},\cd,
-\lm_{-j+i+2}-\lm_{-j+i+3}-\cd -\lm_{-j+n+1}),\\
&& \q\q\q\q\q\q -\lm_{-j+i+2}-\lm_{-j+i+3}-\cd -\lm_{-j+n+2})
\end{eqnarray*}

Now, we set 
\begin{eqnarray*}
X:&=&{\rm max}(0,-\lm_{-j+i+2},-\lm_{-j+i+2}-\lm_{-j+i+3},\cd,
-\lm_{-j+i+2}-\lm_{-j+i+3}-\cd -\lm_{-j+n+1}) \geq 0,\\
Y:&=&-\lm_{-j+i+2}-\lm_{-j+i+3}-\cd -\lm_{-j+n+2}.
\end{eqnarray*}
Then,
$$D_{-j;i}=X+(-\lm_{-j+i+1}+{\rm max}(X,Y))_+-{\rm max}(X,Y)$$

\vskip5pt
(i-a) The case ${\rm max}(X,Y)=X $.

\noindent
We have
\begin{eqnarray*}
D_{-j;i}&=&(-\lm_{-j+i+1}+X)_+\\
&=&(-\lm_{-j+i+1}+{\rm max}(0,-\lm_{-j+i+2},-\lm_{-j+i+2}-\lm_{-j+i+3},\cd,
-\lm_{-j+i+2}-\lm_{-j+i+3}-\cd -\lm_{-j+n+1}))_+\\
&=&{\rm max}(0,-\lm_{-j+i+1},-\lm_{-j+i+2}-\lm_{-j+i+2},\cd,
-\lm_{-j+i+1}-\lm_{-j+i+2}-\cd -\lm_{-j+n+1})\\
&=&C_{-j;i}.
\end{eqnarray*}
In this case, the lemma is true.

\vskip5pt
(i-b) The case ${\rm max}(X,Y)=Y $.

\noindent
In this case, we have
\begin{eqnarray*}
%X&=&0 \text{ or } -\lm_{-j+i+2}-\lm_{-j+i+3}-\cd -\lm_{-j+k+1}
%\text{ for some }k\,(i+1\leq k \leq n),\\
%X-Y&=&\lm_{-j+k+1}+\lm_{-j+k+2}+\cd+\lm_{-j+n+2}<0,\\
D_{-j;i}&=&X-Y+(-\lm_{-j+i+1}+Y)_+.
%{\rm max}(-\lm_{-j+i+1}-\lm_{-j+i+2}-\cd -\lm_{-j+k+1},)
\end{eqnarray*}
First, we consider the case $X=0$. If $D_{-j;i}>0$, then
$-\lm_{-j+i+1}>0$. By the definition of $\lm$, 
we have $\lm_{-j+i+1},\lm_{-j+i+2},\cd, \lm_{-j+n+1}<0$.
This contradicts $X=0$.
Then, we consider the only case $D_{-j;i}\leq 0$. 
In this case, we have $-\lm_{-j+i+1} \leq 0$.
It follows from $X=0$ that
\begin{eqnarray*}
&&\begin{cases}
 -\lm_{-j+i+1} \leq 0,\\
 -\lm_{-j+i+1}-\lm_{-j+i+2} \leq 0,\\
 \cd \\
 -\lm_{-j+i+1}-\lm_{-j+i+2}\cd-\lm_{-j+n+2} \leq 0.
\end{cases}
\end{eqnarray*}
This shows $ C_{-j;i}=0$.
Therefore, since $X\geq 0$, we consider the case $X>0$.

There exist $k$ such that
\begin{eqnarray*}
X&=&-\lm_{-j+i+2}-\lm_{-j+i+3}-\cd-\lm_{-j+k+1}>0.
\end{eqnarray*}
Then we have
\begin{eqnarray*}
X-Y&=&\lm_{-j+k+2}+\lm_{-j+k+3}+\cd+\lm_{-j+n+2}<0,\\
D_{-j;i}&=&X-Y+(-\lm_{-j+i+1}+Y)_+\\
        &=&{\rm max}(X-Y, -\lm_{-j+i+1}+X),\\
C_{-j;i}&=&{\rm max}(0, -\lm_{-j+i+1}+{\rm max}(
-\lm_{-j+i+2},-\lm_{-j+i+2}-\lm_{-j+i+3},\cd-\lm_{-j+i+2}-\cd-\lm_{-j+i+1}))\\
&=&{\rm max}(0, -\lm_{-j+i+1}+X).
\end{eqnarray*}
If $D_{-j;i} \leq 0$, then $-\lm_{-j+i+1}+X \leq 0$. This shows
$C_{-j;i}=0$. If $D_{-j;i} > 0$,  since $X-Y<0$ we have 
$-\lm_{-j+i+1}+X > 0$, then   
$D_{-j;i}=-\lm_{-j+i+1}+X=C_{-j;i}$.
Therefore, if $D_{-j+1;i}>0$, we have shown the lemma.
\vskip5pt
(ii) The case $D_{-j+1;i} \leq 0 $.

\noindent
In this case, $C_{-j+1;i}=0$. 
Then, we have
\begin{eqnarray*}
&&\begin{cases}
 -\lm_{-j+i+2} \leq 0,\\
 -\lm_{-j+i+2}-\lm_{-j+i+3} \leq 0,\\
 \cd \\
 -\lm_{-j+i+2}-\lm_{-j+i+3}\cd-\lm_{-j+n+2} \leq 0.
\end{cases}
\end{eqnarray*}
This shows, by the definition of $\lm$, $\lm_{-j+i+1}\geq 0$ and
$-\lm_{-j+i+1}+X \leq 0$.
Since $0=C_{-j+1;i} \geq C_{-j;i+1} \geq 0$ we have $C_{-j;i+1}=0$.
Then, 
\begin{align*}
D_{-j;i}&=C_{-j;i+1}+C_{-j+1;i-1}-C_{-j+1;i}+D_{-j+1;i}\\
&=(-{\lm_{-j+i+1}}+{\rm max}(0,-\lm_{-j+i+2},\cd,
-\lm_{-j+i+2}-\lm_{-j+i+3}-\cd -\lm_{-j+n+2}))_++D_{-j+1;i}\\
&=(-\lm_{-j+i+1}+X)_++D_{-j+1;i}\\
&=D_{-j+1;i} \leq 0.
\end{align*}
On the other hand, 
\begin{eqnarray*}
C_{-j;i}&=&{\rm max}(0, -\lm_{-j+i+1}+{\rm max}(
-\lm_{-j+i+2},-\lm_{-j+i+2}-\lm_{-j+i+3},\cd-\lm_{-j+i+2}-\cd-\lm_{-j+i+1}))\\
&=&{\rm max}(0, -\lm_{-j+i+1}+X)=0.
\end{eqnarray*}

This completes the proof of Lemma $\ref{lem-d}$. 
\qed

\vskip5pt
Let $v_0=(\cd,0,0,t_\lm,x_{-1;n},x_{-1;n-1},\cd,x_{-j;i},\cd)$ 
be a highest weight vector in $\Sigma_\io[\lm]\cap\Sigma'_\io[\lm]$
satisfying :
\begin{equation}
x_{-j;i}+C_{-j;i}\geq0.
\label{xc}
\end{equation}
Note that the linear function $x_{-j;i}+C_{-j;i}$ in the 
left hand-side of the inequality is a generator of $\Xi'_\io[\lm]$.
Thus, any vector in $\Sigma_\io[\lm]\cap\Sigma'_\io[\lm] $
satisfys the inequality (\ref{xc}).

Now, we shall show that $v_0$ is uniquely determined and 
coincides with $v_0$ by the induction on the index $(-j;i)$,
where we consider the lexicographic order :
$(-j;i)<(-j';i')$ if $j<j'$, or $j=j'$ and $i>i'$,{\it i.e.},
\[
(-1;n)<(-1;n-1)<\cd<(-1;1) <(-2;n)<\cd<\cd (-j;i+1)
<(-j;i)<(-j;i-1)<\cd.
\]
Recall the condition that $x_0$ is a highest weight vector :
\[
 \sigma_{-j;i}(v_0)\leq 0 {\rm \,\,for \,\,any }(-j;i)\,\,
(j\geq1,1\leq i\leq n).
\]
For $(-j;i)=(-n;1)$, we have 
$\sigma_{-1;n}(v_0)=x_{-1;n}-\lm_n\leq0$.
By (\ref{xc}), we also have 
$x_{-1;n}+C_{-1;n}=x_{-1;n}-\lm_n\geq0$.
Those imply 
\[
 x_{-1;n}=\lm_n=-C_{-1;n}.
\]
Here note that by the assumption on $\lm$, we have $\lm_n\leq0$.
Asuume that for any $(-j';i')< (-j;i)$, 
\begin{equation}
x_{-j';i'}=-C_{-j';i'}.
\label{ass}
\end{equation}
Let us determine $x_{-j;i}$.
By this assumption (\ref{ass}), we have
\begin{equation}
\sigma_{-j;i}(v_0)=x_{-j;i}+D_{-j;i}.
\label{xd}
\end{equation}
If $D_{-j;i}\leq 0$, by Lemma \ref{lem-d}(i) we have
$C_{-j;i}=0$. Then we have 
\[
0\leq x_{-j;i}+C_{-j;i}=x_{-j;i}.
\]
On the other hand, 
since $v_0$ is an element in $\Sigma_\io[\lm]\cap\Sigma'_\io[\lm]$,
we have $x_{-j;i}\leq0$. Those impliy $x_{-j;i}=0=-C_{-j;i}$.

If $D_{-j;i}> 0$, by Lemma \ref{lem-d}(ii)
 we have $D_{-j;i}=C_{-j;i}$. 
Since $v_0$ is a highest weight vector, we have
\[
0\geq\sigma_{-j;i}(v_0)=x_{-j;i}+D_{-j;i}=x_{-j;i}+C_{-j;i}.
\]
By the condition (\ref{xc}), we have 
\[
 x_{-j;i}+C_{-j;i}\geq0.
\]
Thus, we obtain $x_{-j;i}=-C_{-j;i}.$
Now, we know that $v_0$ is the unique highest weight 
vector in $\Sigma_\io[\lm]\cap\Sigma'_\io[\lm]$ satisfying 
$x_{-j;i}+C_{-j;i}\geq0$. Since $B_0(\lm)$ contains 
the unique highest weight vector (\cite{K3}), 
$v_0$ must be the unique highest weight vector in $B_0(\lm)$,
which implies that 
$B_0(\lm)=\Sigma_\io[\lm]\cap \Sigma'_\io[\lm]$.
\qed

%%%%%%%%% Section 6   %%%%%%%%%%
\section{Polyhedral Realization of $B(\uq a_{\lm})$ of Type $A^{(1)}_1$}
\setcounter{equation}{0}
\renewcommand{\theequation}{\thesection.\arabic{equation}}
%%%%%%%%%%%%%%%%%%%%
In this section, we consider the case $\ge$ is of type $A^{(1)}_1$.
We fix a positive level integral weight 
$\lm=\lm_1\Lm_1+\lm_2\Lm_2$ $(\lm_1\in \ZZ_{ >0}, \lm_2\in
\ZZ_{\leq0},\,\lm_1+\lm_2>0)$.
We define $C_{-k}$ $( k \in \ZZ_{>0})$ as follows:
$$ C_{-k} := (-(k-1)\lm_1 -k\lm_2 )_{+}.$$
\begin{thm1}
\label{A_1}
Let $\io = (\cd,2,1,2,1,t_{\lm},2,1,2,1,\cd)$ be an infinite 
sequence and let $B_0(\lm)$ be the connected component of 
${\rm Im} \, (\Psi^{(\lm)}_{\io })$ containing $\vec 0 := 
 (\cd,0,0,t_{\lm},0,0,\cd)$. Set
\begin{align*}
&\Xi_{\io}[\lm]:=\{\bar{S}_{j_l}\cd \bar{S}_{j_1}(x_{j_0})\,
  :\,l\geq0,\,j_0,\cd,j_l\ne 1\}\\
  &\q\q\;\; \cup\{\bar{S}_{{-j}_k}\cd \bar{S}_{{-j}_1}(-x_{{-j}_0})\,
  :\,k\geq0,\,j_0,\cd,j_k\ne1\},\\
&\Sigma_{\io}[\lm]
 :=\{\vec x\in \ZZ^{\ify}_{\io}[\lm](\subset \QQ^{\ify})\,:\,
\vp(\vec x)\geq 0\,\,{\rm for \,\,any }\,\,\vp\in \Xi_{\io}[\lm]\},\\
&\Xi^{'}_{\io}[\lm]:=
\{\bar{S}_{-{j_l}}\cd \bar{S}_{-{j_1}}({x_{-k}}+C_{-k})\,
:\,l\geq0,\,k\geq1 ,\ j_1,\cd,j_l \geq 1\},\\
&\Sigma^{'}_{\io}[\lm]:=
    \{\vec x\in \ZZ^{\ify}_{\io}[\lm](\subset \QQ^{\ify})\,:\,
    \vp(\vec x)\geq 0\,\,{\rm for \,\,any }\,\,\vp\in
 \Xi^{'}_{\io}[\lm]\}.
\end{align*}
Then,
  \begin{enumerate}
   \item
    $B_0(\lm) = \Sigma_{\io}[\lm]\cap\Sigma{'}_{\io}[\lm].$
   \item
Let $v_{\lm}$ be the highest weight vector of
 $B_0(\lm)$. Then we have
 $$v_{\lm} = (\cd,0,0,t_{\lm},-C_{-1},-C_{ -2},\cd,-C_{-k},\cd).$$
   \end{enumerate}
\end{thm1}
{\sl Proof.}\,\,
Since $\Xi_{\io}[\lm]$ is closed by $\bar{S}_{k}$'s,
 by Lemma 4.3 it has a crystal structure. We will show that 
$\Sigma_{\io}[\lm]\cap\Sigma{'}_{\io}[\lm]$ contains $\vec 0$ 
and has the unique highest weight vector. 
First, we will show that
\begin{eqnarray*}
  x_0 =(\cd,0,0,t_{\lm},-C_{-1},
    -C_{ -2},\cd,-C_{-k},\cd)\,\,
   (1\le k)
\end{eqnarray*}
is one of the highest weight vectors of 
$\ZZ^\ify_\io[\lm]$.
For $k\geq1$, we recall the definition of $\sigma_{-k}$:
$$\sigma_{-k}(\vec x)= 
-\lan h_{i_{-k}},\lm\ran+x_{-k}+\sum_{j>-k}\lan h_{i_{-k}},\al_{i_j}\ran x_j.$$
\vskip5pt
Let us show $\sigma_{-k}(x_0)\leq0$ for $k\geq1$ 
by the induction on $k$. ( In the case $k<0$, trivially 
$\sigma_{-k}(x_0) = 0$. )

\vskip5pt
(i) The case $k=1$.

\noindent
we have $\sigma_{-1}(x_0) = -(-\lm_2)_+ - \lm_2 = 0$.
\vskip5pt

(ii) The case $k=2$.

\noindent
By $(\ref{zeroika})$, we have 
\begin{eqnarray*}
  \sigma_{-2}(x_0) &=& -(-\lm_1-2\lm_2)_+ + 2(-\lm_2)_+ - \lm_1\\
            &=& -(-\lm_1-2\lm_2)_+ - 2\lm_2 - \lm_1\,\,(by\,\,\lm_2\leq0)\\
            &\leq0.&
\end{eqnarray*}
\vskip5pt

(iii) The case $k>2$.

\noindent
We assume $\sigma_{-k+2}(x_0)\leq0$. The following fact is trivial
 by $\lm_1+\lm_2>0$:
\begin{eqnarray}
   -(k-1)\lm_1 - k\lm_2 \leq 0 &\Longrightarrow& -k\lm_1 -(k+1)\lm_2 \leq 0.
\label{trivial}
\end{eqnarray}
By the definition of $\sigma_{-k}(x_0)$, we heve
$$\sigma_{-k}(x_0)=-(-(k-1)\lm_1-k\lm_2)_++2(-(k-2)\lm_1-(k-1)\lm_2)_+
-(-(k-3)\lm_1-(k-2)\lm_2)_++ \sigma_{-k+2}(x_0).$$
Here, we set 
$$X:=-(k-1)\lm_1-k\lm_2,\, Y:=-(k-2)\lm_1-(k-1)\lm_2,\, 
Z:=-(k-3)\lm_1-(k-2)\lm_2.$$
Then we have
 $\sigma_{-k}(x_0)=-(X)_++2(Y)_+-(Z)_++\sigma_{-k+2}(x_0)$.
\vskip5pt
By $(\ref{trivial})$, if $Z\leq 0$, then $Y \leq 0$, $X\leq 0$ and 
if $Y \leq 0$, then $X\leq 0$. It is sufficient to show following 
four cases:

(iii-a) $Z\leq 0$.\,\,(iii-b) $Z>0$, $Y\leq0$.\,\, 
(iii-c) $Z>0$, $Y>0$, $X\leq0$.\,\,
(iii-d)  $Z>0$, $Y>0$, $X>0$.

\vskip5pt
(iii-a) The case $Z\leq 0$.

\noindent
We have
$\sigma_{-k}(x_0)=\sigma_{-k+2}(x_0)\leq0$.
\vskip5pt

(iii-b) The case $Z>0$, $Y\leq0$.

\noindent
By definition of $\sigma_{-k}(x_0)$,
we have $\sigma_{-k}(x_0)=-Z+\sigma_{-k+2}(x_0)\leq0$.
\vskip5pt

(iii-c) The case $Z>0$, $Y>0$, $X\leq0$.
\begin{align*}
\sigma_{-k}(x_0)&=2(-(k-2)\lm_1-(k-1)\lm_2)_+-(-(k-3)\lm_1-(k-2)\lm_2)_+
+\sigma_{-k+2}(x_0)\\
&=X+\sigma_{-k+2}(x_0)\leq0.
\end{align*}

(iii-d) The case $Z>0$, $Y>0$, $X>0$.

$\sigma_{-k}(x_0)= -X +2Y -Z + \sigma_{-k+2}(x_0)=\sigma_{-k+2}(x_0)\leq 0$.
\vskip5pt

Next, we will show that $\Sigma_{\io}[\lm]\cap\Sigma{'}_{\io}[\lm]$
contains $\vec 0$.
By the similar way to 
the proof of Theorem$\ref{thm1}$, we will show $\vec 0 \in
\Sigma{'}_{\io}[\lm]$.
For the purpose, we shall calculate the constant term of 
$\bar{S}_{-{j_k}}\cd\bar{S}_{-{j_1}}({x_{-k}})$. 

We set
\begin{equation*}
\vp^{(l)}_{-k}(x):=
\begin{cases}
\bar{S}_{-k+l-1}\cd\bar{S}_{-k+1}\bar{S}_{-k}(x_{-k})&(l\leq k),\\
\bar{S}_{-k+l}\cd\bar{S}_{1}\bar{S}_{-1}\bar{S}_{-2}
\cd\bar{S}_{-k+1}\bar{S}_{-k}(x_{-k})&(l>k).
\end{cases}
\end{equation*}
By the similar argument in $\cite{NZ}$ Lemma4.2, 
we obtain the explicit form of
$\vp^{(l)}_{-k}(x)$ up to constant term as follows:

\begin{lem1}
\begin{eqnarray*}
&&\Xi^{'}_{\io}[\lm]=\{\vp^{(l)}_{-k}(x)+C_{-k}\,|\,k \geq
1,\,l\geq1\,\},\\
&&\vp^{(l)}_{-k}(x)-\vp^{(l)}_{-k}(0)=
(l+1)x_{l-k+\theta(l-k)}-lx_{l-k+1+\theta(l-k+1)}
\end{eqnarray*}
where
\[
\theta(x):=
\begin{cases}
1 &\text{if $x\geq0$,}\\
0 &\text{if $x<0$.}
\end{cases}
\]
\end{lem1} 

Now, we calculate the constant term in $\vp^{(l)}_{-k}$.
For $l \leq k-2$, we know that $\vp^{(l)}_{-k}$ has no costant
term by its definition. And we have
$$\vp^{(k-2)}_{-k}(x)=(k-1)x_{-2}-(k-2)x_{-1}.$$
By direct calculations, we obtain
\begin{eqnarray*}
\vp^{(k-1)}_{-k}&=&\bar{S}_{-2}(\vp^{(k-2)}_{-k})\\
                &=&kx_{-1}-(k-1)x_{1}+(k-1)\lm_{1},\\
\vp^{(k)}_{-k}\q  &=&\bar{S}_{-1}\bar{S}_{-2}(\vp^{(k-2)}_{-k})\\
                &=&(k+1)x_1-kx_2+(k-1)\lm_1+k\lm_2. 
\end{eqnarray*}

For $k\geq1$, since $\bar{S}_{k}$ does not produce non-trivial
constant term, we have
$$\vp^{(k)}_{-k}(0)=\vp^{(l)}_{-k}(0)\,\,(l>k).$$
Hence, we obtain the constant term of $\vp^{(l)}_{-k}(x)$:
\begin{equation}
\vp^{(l)}_{-k}(0)=
\begin{cases} 
(k-1)\lm_1+k\lm_2 &\,\,(l\geq k),\\
(k-1)\lm_1 &\,\,(l=k-1),\\
0 &\,\,(l\leq k-2).
\end{cases}
\end{equation}
By Lemma $\ref{r1-n}$, we have $C_{-k}+\vp^{(l)}_{-k}(0)\geq 0$.
This shows that constant terms in all elements in $\Xi^{'}_{\io}[\lm]$
are non-negative and then $\vec 0$ is contained in 
$\Sigma_{\io}[\lm]\cap\Sigma^{'}_{\io}[\lm]$.
Therefore, we have 
$B_0(\lm)\subset \Sigma_{\io}[\lm]\cap\Sigma{'}_{\io}[\lm]$.

Finally, we will show that 
$\Sigma_{\io}[\lm]\cap\Sigma{'}_{\io}[\lm]$ has 
the unique highest weight vector $v_{\lm}$.

We define
\begin{equation*}
 D_{-k}:= \sigma_{-k}(x_0)+C_{-k}.
\end{equation*}

Note that for $k \geq 3$, 
 \begin{align}
  D_{-k} &= 2(-(k-2)\lm_1-(k-1)\lm_2)_+ -(-(k-3)\lm_1-(k-2)\lm_2)_+
             + \sigma_{-k+2}(x_0) \nn \\
         &=2(Y)_+ -(Z)_++ \sigma_{-k+2}(x_0),
\label{defd}\\
  C_{-k} &= (-(k-1)\lm_1-k\lm_2)_+ \nn\\
         &= (X)_+.
  \label{defc}
 \end{align}

We need the following lemma:

 \begin{lem1}
  \begin{enumerate}
   \item If $D_{-k}\leq 0$, then $C_{-k}=0$.
   \item If $D_{-k}> 0$, then we have $D_{-k}=(-(k-1)\lm_1-k\lm_2)_+=C_{-k}$.
  \end{enumerate}
 \label{lemd}
 \end{lem1}
{\sl Proof.} We shall show the lemma by the induction on $k$.

\nd
(I) The case $k=1$. 

\noindent
By the definition of $D_{-1}$, we have $D_{-1} = -\lm_2 > 0$.
This shows $C_{-1}=(-\lm_2)_+ = -\lm_2$.

\nd
(II) The case $k=2$.

\noindent
By the definition of $D_{-2}$ and $C_{-2}$, we have
\[
   D_{-2} = 2(-\lm_2)_+ -\lm_1 = -2\lm_2 -\lm_1,\q
   C_{-2} = (-\lm_1-2\lm_2)_+ = (D_{-2})_+.
\]
Therefore, we get if $D_{-2}\leq 0$, then $C_{-2}=0$.
If $ D_{-2}> 0$,  then $ D_{-2} =-\lm_1-2\lm_2= C_{-2}$.

\nd
(III) The case $k\geq 3$.

\nd
(i) The case $D_{-k+2}>0$.

\noindent
By the induction hypothesis, $\sigma_{-k+2}(x_0)=D_{-k+2}-C_{-k+2}=0$.
By (\ref{defd}) and (\ref{defc})
\[
D_{-k}= 2(Y)_+-(Z)_+,\qq C_{-k}=(X)_+.
\]
By $(\ref{trivial})$, if $Z\leq 0$, then $Y\leq 0$ and $X\leq 0 $. 
In this case, we have
$D_{-k}=0$ and $C_{-k}=0$. 
If $Y\leq 0$, then $X\leq 0 $. In this case, $D_{-k}=-(Z)_+ \leq 0$
and $C_{-k}=0$.
Then, it is sufficient to show following two cases:

(i-a) $Z>0$, $Y>0$ and $X \leq 0$. 

(i-b) $Z>0$, $Y>0$ and $X> 0$.

\nd
(i-a) The case $Z>0$, $Y>0$ and $X \leq 0$.

\noindent
We have
\begin{eqnarray*}
D_{-k}&=&2(-(k-2)\lm_1-(k-1)\lm_2)-(-(k-3)\lm_1-(k-2)\lm_2)\\
      &=&-(k-1)\lm_1-k\lm_2=X \leq 0,\\
C_{-k}&=&(X)_+=0.
\end{eqnarray*}
This show that the lemma is true.

\nd
(i-b) The case $Z>0$, $Y>0$ and $X> 0$.

\noindent
We have
\begin{eqnarray*}
D_{-k}&=&2(-(k-2)\lm_1-(k-1)\lm_2)-(-(k-3)\lm_1-(k-2)\lm_2)\\
      &=&-(k-1)\lm_1-k\lm_2=X>0,\\
C_{-k}&=&(X)_+=X =D_{-k}.
\end{eqnarray*}
Therefore, in the case $D_{-k+2}>0$, the lemma is true.

\nd
(ii) The case $D_{-k+2}\leq 0$.

\noindent
By the induction hypothesis, $C_{-k+2}=(-(k-3)\lm_1-(k-2)\lm_2)_+=0$,
and then $Z=-(k-3)\lm_1-(k-2)\lm_2 \leq 0$. By $(\ref{trivial})$, 
$Y$, $X \leq 0$. By (ref{defd}) and (\ref{defc}), we have
\[
   D_{-k} = \sigma_{-k+2}(x_0) \leq 0,\qq
  C_{-k} = 0.
\]
Therefore, we complete the proof.
\qed

Let $v_0:=(\cd,0,0,t_{\lm},x_{-1},x_{-2},\cd,x_{-k},\cd)$ 
be a highest weight vector in 
$\Sigma_{\io}[\lm]\cap\Sigma^{'}_{\io}[\lm]$, which satisfys:
\begin{eqnarray}
x_{-k}+C_{-k}\geq0.
\label{ineq}
\end{eqnarray}
Note that the linear function $x_{-k}+C_{-k}$ is a generator 
of $\Xi^{'}_{\io}[\lm]$, and then
 any vector in
$\Sigma_{\io}[\lm]\cap\Sigma^{'}_{\io}[\lm]$
satisfys the inequality $(\ref{ineq})$.

We shall show that $v_0$ is uniquely determined and coincides with
$v_{\lm}$ by the induction on the index $k$.

By the condition that $v_0$ is a highest weight vector, we have 
that $\sigma_{-k}(v_0)\leq0 {\rm \,\,for \,\,any\,\, } k \geq 1.$
For $k=1$, we have 
$\sigma_{-1}(v_0)=x_{-1}-\lm_2\leq0$.
By (\ref{ineq}), we also have 
$x_{-1}+C_{-1}=x_{-1}-\lm_2\geq0$,
which implies $x_{-1}=\lm_2=-C_{-1}$.
Assume that for any $k'< k$, 
\begin{equation}
x_{-k'}=-C_{-k'}.
\label{ass2}
\end{equation}
Let us determine $x_{-k}$.
By the assumption (\ref{ass2}), we have
\begin{equation}
\sigma_{-k}(v_0)=x_{-k}+D_{-k}.
\label{xd2}
\end{equation}
If $D_{-k}\leq 0$, by Lemma \ref{lemd}(i) we have
$C_{-k}=0$. Then we have 
$0\leq x_{-k}+C_{-k}=x_{-k}$.

\nd
On the other hand, 
since $v_0$ is an element in $\Sigma_\io[\lm]\cap\Sigma'_\io[\lm]$,
we have $x_{-k}\leq0$. We obtain $x_{-k}=0=-C_{-k}$.

If $D_{-k}> 0$, by Lemma \ref{lemd}(ii)
 we have $D_{-k}=C_{-k}$. 
Since $v_0$ is a highest weight vector, we have
\[
0\geq\sigma_{-k}(v_0)=x_{-k}+D_{-k}=x_{-k}+C_{-k}.
\]
On the other hand, by the condition (\ref{ineq}), we have 
$x_{-k}+C_{-k}\geq0$.
Thus, we obtain $x_{-k}=-C_{-k}.$
Now, we know that $v_0$ is the unique highest weight 
vector in $\Sigma_\io[\lm]\cap\Sigma'_\io[\lm]$ satisfying 
$x_{-k}+C_{-k}\geq0$. Since $B_0(\lm)$ contains 
the unique highest weight vector (\cite{K3}), 
$v_0$ must be the unique highest weight vector in $B_0(\lm)$.
Hence, we have $B_0(\lm)=
\Sigma_{\io}[\lm]\cap\Sigma{'}_{\io}[\lm]$.
\qed

\end{document}